\documentclass[review]{elsarticle}
\usepackage[utf8]{inputenc}

\usepackage{lineno,hyperref}

\usepackage{amsmath}
\usepackage{amsfonts} 
\usepackage{longtable}
\usepackage{lscape}
\usepackage{comment}
\usepackage[official]{eurosym}
\usepackage[dvipsnames]{xcolor}
\modulolinenumbers[5]
\definecolor{dkred}{rgb}{0.8,0,0}
\definecolor{blue}{rgb}{0,0,1}

\usepackage{circuitikz}

\journal{Applied Energy}

\begin{document}

\begin{frontmatter}

\title{Ramping up the hydrogen sector: An energy system modeling framework}

\author[mymainaddress]{T. Klatzer\corref{mycorrespondingauthor}}
\cortext[mycorrespondingauthor]{Corresponding author:}
\ead[url]{thomas.klatzer@tugraz.at}
\author[mymainaddress]{U. Bachhiesl}
\author[mymainaddress]{S. Wogrin}
\author[mysecondaryaddress]{A. Tomasgard}

\address[mymainaddress]{Institute of Electricity Economics and Energy Innovation (IEE), Graz University of Technology, Inffeldgasse 18, Graz, Austria}
\address[mysecondaryaddress]{Department of Industrial Economics and Technology Management (I\O T), The Norwegian University of Science and Technology, Sentralbygg 1, Gl\o shaugen, Trondheim, Norway}

\begin{abstract}
With the transition towards a decarbonized society, energy system integration is becoming ever more essential. In this transition, the energy vector hydrogen is expected to play a key role as it can be produced from (renewable) power and utilized in a plethora of applications and processes across sectors. To date, however, there is no infrastructure for the production, storage, and transport of renewable hydrogen, nor is there a demand for it on a larger scale. In order to link production and demand sites, it is planned to re-purpose and expand the existing European gas pipeline network in the future. 
During the early stages of ramping up the hydrogen sector (2020s and early 2030s), however, blending natural gas with hydrogen for joint pipeline transmission has been suggested. 
Against this background, this paper studies hydrogen blending from a modeling perspective, both in terms of the implications of considering (or omitting) technical modeling details and in terms of the potential impact on the ramp-up of the hydrogen sector.
To this end, we present a highly modular and flexible integrated sector-coupled energy system optimization model of the power, natural gas, and hydrogen sectors with a novel gas flow formulation for modeling blending in the context of steady-state gas flows.
A stylized case study illustrates that hydrogen blending has the potential to initiate and to facilitate the ramp-up of the hydrogen sector, while omitting the technical realities of gas flows -- particularly in the context of blending -- can result in suboptimal expansion planning not only in the hydrogen, but also in the power sector, as well as in an operationally infeasible system.
\end{abstract}


\end{frontmatter}

\section{Introduction}
\label{sec:Intro}

\subsection{Motivation}
\label{sec:Motivation}
The transition to carbon-neutral societies requires both societal changes and the transformation of energy systems. The European Union (EU) is on the frontier of this process and has set the goal of achieving climate neutrality by 2050 \cite{GreenDeal2019}. 
In light of Russia's invasion of Ukraine, however, energy security and affordability have come into focus, at least in the short term. For the medium term, the interim targets for 2030 according to Fit for 55 \cite{Fitfor55} are tightened under the RePowerEU plan \cite{RePowerEUPlan2022}. In addition to the massive expansion of renewable energy capacities, the dependence on natural gas is to be reduced and the ramp-up of the hydrogen sector is to be accelerated. The goal is to produce 10~Mt of renewable\footnote{The legal framework defining when hydrogen qualifies as renewable is still missing. To date, only a draft delegated regulation \cite{DefRenewableH2Draft2022} exists.} hydrogen in the EU alongside 6~Mt of renewable hydrogen (and 4~Mt of ammonia) imports by 2030. To implement the RePowerEU plan \cite{ImplementingREPowerEU2022}, 65~GW of electrolyzers powered by 41~GW wind and 62~GW solar are to be built, which involves highly complex techno-economic issues that could be described as a chicken-and-egg problem.

For example, it is unclear where hydrogen production infrastructure should ideally be sited -- in proximity to renewables or in proximity to (future) hydrogen demand?
Especially the latter is still subject to a lot of uncertainty. This is since the energy vector hydrogen offers the potential for a plethora of applications, ranging from the (seasonal) storage of electricity from variable renewable energy (VRE) sources, to its use as a feedstock, e.g. for the direct reduction of iron and for e-fuels, to its energetic use as a substitute for natural gas.
In addition to production, the capability to transport hydrogen will have a significant impact on how quickly the sector will evolve.
Particularly during the early stages of ramping up the hydrogen sector (2020s to early 2030s), utilizing the well-established European pipeline transmission system to transport hydrogen by blending it with natural gas has been found a promising mode of transportation \cite{Jens2021}. Furthermore, blending is already permitted to date, e.g. up to 10\% in Austria \cite{GasConnect2021}, and potentially up to 20\% in the United Kingdom in the near future \cite{UKH2Strategy2021}. Against this background, assessing hydrogen blending in the context of expansion planning is highly relevant.

Motivated by the issues stated above, the research questions studied in this paper are: 
\begin{itemize}
    \item How to model natural gas and hydrogen blending in integrated  sector-coupled energy system models (ESMs) and what are the implications of the technical modeling details on expansion planning?
    \item Leveraging existing natural gas infrastructure for blending and policies, how could the hydrogen sector ramp up in a holistic and cost-effective way, especially during its early stages when there is no dedicated hydrogen demand yet?
\end{itemize}

To address these research questions, in this paper we propose a novel gas flow framework for detailed modeling of natural gas and hydrogen blending based on steady-state gas flows. To test the proposed framework, we formulate an integrated sector-coupled ESM for expansion planning in the power, natural gas, and hydrogen sectors, with particular emphasis on the technical modeling aspects of the natural gas and hydrogen sectors (production, transmission, and the demand side).
In our case studies we illustrate that a high level of technical detail, e.g. the formulation of gas flows and blending, is valuable for ESMs as omitting these can significantly impact operational feasibility of the planned energy system.
Moreover, in integrated sector-coupled ESMs capturing techno-economic details can be a decisive factor in determining to which extent and in which sectors a technology should be primarily deployed in the future. 

In the following section, we conduct a literature review on state-of-the-art ESMs of the power and natural gas sectors and study whether the research questions posed can be addressed with existing ESMs.

\subsection{Literature review}
\label{sec:Literature}
This literature review addresses existing ESMs of the power and natural gas sectors, and their approach for modeling expansion of hydrogen infrastructure during the early stages of ramping up the hydrogen sector.
In particular, we focus on the formulation of gas flows, which is key for modeling detailed natural gas and hydrogen blending.
Since this only represents a specific segment of ESMs, the interested reader is also referred to other more comprehensive literature reviews, e.g. \cite{Pfenninger2014,Fodstad2022}. 
For example, Pfenninger et al. \cite{Pfenninger2014} examine different modes of ESMs in terms of their capacity to address the energy challenges of the twenty-first century and find that traditional ESMs might not be capable of deriving feasible solutions for increasingly complex systems under the decarbonization paradigm. 
Fodstad et al. \cite{Fodstad2022} point out that ESMs that consider power and natural gas systems typically focus on either reliability and security of supply aspects (in the short term) or system expansion (in the long term). However, a distinct trade-off between representation of technical detail and temporal resolution remains.
Furthermore, they identify hydrogen value chains as concrete modeling frontiers.

Let us now take a closer look at existing ESMs and discuss them with a focus on power, natural gas, and hydrogen modeling.
METIS \cite{METIS2018} is a linear program (LP) energy modeling tool for policy making that considers power, natural gas, hydrogen, and bio-methane to decarbonize Europe's energy system \cite{Arduin2022}. METIS represents all energy flows exclusively as transport problems (TPs) - at least in its standard version.
PyPSA-Eur-Sec \cite{PyPSA-Eur-Sec2022} is a powerful open-source ESM extension of PyPSA-Eur \cite{Horsch2018}, a European power system expansion planning model comprising the transport, heating, and industry sectors. The power flow is formulated as linear optimal power flow \cite{Horsch2018a}, while the hydrogen network is a TP (greenfield approach) that currently does not include the option of gas blending. For future releases it is planned to include the existing natural gas network as TP. PyPSA-Eur-Sec is an LP and therefore does not allow for integer variables, e.g. for modeling of unit commitment (UC) decisions. However, UC is likely to be relevant for future ESMs, e.g. for operation of hydrogen-fired power plants\footnote{Industry is working on a classification scheme on hydrogen-readiness for gas power plants~\cite{EUT2021}.}. The model has been applied for case studies, e.g. in \cite{Victoria2022}.
PRIMES \cite{PRIMES2018} is an extensive long-term ESM based on 5-year periods.
PRIMES provides several sub-models. This includes a 'power and steam generation and supply model' with a direct current optimal power flow (DC-OPF) (on country level), a 'gas supply model' based on a single node (SN) representation per country connected by a TP, a 'new fuels and storage model' to study hydrogen, Power-to-X, synthetic fuels etc., and many more. The latter is based on hourly resolution and allows for distinct natural gas and hydrogen transmission as well as blending. However, PRIMES is a commercial model and the mathematical formulation for the respective sub-model is not publicly available\footnote{We assume that both distinct natural gas and hydrogen transmission as well as blending are also based on TPs (see 'gas supply model').}.

While the above frameworks comprise of continental scales, there are several models that focus on the national level, e.g. \cite{Husarek2021,Gils2021}. 
Husarek et al. \cite{Husarek2021} extend a holistic, dynamic ESM (based on hourly resolution and the NUTS 2 level) to investigate hydrogen supply scenarios for Germany up to 2050. Similarly, Gils et al. \cite{Gils2021} examine the expansion of hydrogen infrastructure in Germany until 2050. The applied REMix model includes a DC-OPF formulation, however, their case study is focused on ten German regions as well as several European countries (SN representation) where hydrogen infrastructure is aggregated per region and the gas transmission system is represented by a TP.

All of the above models consider widely simplified versions of the gas transmission system and disregard physical realities as they are large-scale models.
At the other end of the spectrum, several papers include detailed representations of natural gas networks based on the steady-state gas flow equation, which comprises the pressure drop along a pipeline as well as the direction of gas flows and can either be linearized \cite{Zhang2018,Correa-Posada2015} or applied in a non-linear framework, e.g. \cite{Zeng2017,Tao2021,Qiu2015}. This formulation can be further extended by considering the storage capability of a pipeline, referred to as linepack, e.g. \cite{Correa-Posada2015,Zeng2017,Qiu2015}, which presents a valuable flexibility option for the gas system.
However, the increased level of technical detail comes at the expense of computational complexity and eventually intractability, especially in case of non-linear formulations \cite{Boyd2009}.

From this review (and our proprietary extensive literature review \cite{Klatzer2022}), we note that hydrogen transmission is either based on dedicated hydrogen pipelines modeled as TPs \cite{Victoria2022,Bodal2020,He2021} or based on road transport \cite{Husarek2021,He2021,Li2020}.
Natural gas and hydrogen blending, on the other hand, is not considered in existing ESMs. Hence, there is a lack of model formulations for blending per se, e.g. based on steady-state gas flows, as well as studies on blending and its impact on expansion planning, especially with regard to ramping up the hydrogen sector.

\subsection{Proposed ESM and original contributions}
\label{sec:Contributions}
In this work, we propose an integrated sector-coupled ESM of the power, natural gas, and hydrogen sectors based on the objective of minimizing total system costs (expansion and operation). The presented ESM is formulated as a deterministic mixed-integer linear program (MILP) and is highly modular and flexible, both in terms of modeling features and temporal structure -- all within a unified framework.
This includes, inter alia, the detailed representation of the power \cite{Wogrin2022-LEGO}, natural gas, and hydrogen sectors and their interconnections, e.g. via electrolyzer units, fuel cells, steam-methane reforming units, gas-fired power plants, etc. Gas transport is modeled as a high-pressure gas transmission system in which gas flows can be represented by a (modified) transport problem or as steady-state gas flows. 
As an original modeling contribution, we model natural gas and hydrogen blending for collective pipeline transport based on steady-state gas flows or what we call blending pressure problem. To the best of our knowledge, the high level of technical detail associated with this -- and blending in particular -- is a novelty in the context of ESMs, allowing for more accurate infrastructure planning not only in the gas sector but also in the power sector, especially during the early stages of establishing hydrogen infrastructure.

The remainder of this paper is structured as follows: In section \ref{sec:MathForm} we provide the mathematical framework for the integrated sector-coupled power, natural gas, and hydrogen optimization model. To address the research questions posed, in section \ref{sec:CaseStudies} we conduct two comprehensive case studies based on a stylized integrated power and gas system.
First, we assess the impact of the gas flow formulation -- and blending in particular -- on planning results in the integrated sector-coupled energy system and then we investigate how the hydrogen sector could ramp up in a holistic and cost-effective way by utilizing existing infrastructure. Finally, section~\ref{sec:Conclusions} concludes the paper.

\allowdisplaybreaks
\section{Mathematical formulation of the integrated ESM}
\label{sec:MathForm}
This section provides the mathematical formulation of the integrated power, natural gas, and hydrogen ESM. The formulation represents a significant extension of the Low-carbon Expansion Generation Optimization (LEGO) open-source model \cite{Wogrin2022-LEGO} available on \href{https://github.com/IEE-TUGraz/LEGO}{GitHub}.LEGO is a power system model that includes electrolyzer units, but which does not model detailed gas infrastructure. In this paper, the LEGO model is extended to an integrated sector-coupled ESM. Novelties include:
\begin{itemize}
    \item Hydrogen units: steam-methane reforming, fuel cell, and storage units
    \item Natural gas units: gas wells and storage units
    \item Formulation of pipeline gas flows: as blending transport problem, based on the piecewise linearized steady-state gas flow equation
    \item Hydrogen blending for pipeline gas transmission, to substitute natural gas in the gas sector (including the potential to reduce CO\textsubscript{2} emissions), and for co-firing in gas-fired power plants
\end{itemize}

As a starting point, we briefly introduce the reader to the model's underlying temporal structure in section~\ref{subsec:Time}. Although not a novelty, this is important to get a better understanding of the model and, in particular, gas storage technologies described later on. In the following sections, we present the objective function including bounds (\ref{sec:OF}), elements of the natural gas (\ref{subsec:NG-sector}) and hydrogen (\ref{subsec:H2-sector}) sectors, and the framework governing the gas transmission system (\ref{subsec:GasNetwork}). For the latter and as an original contribution, we consider blending of natural gas and hydrogen in the context of steady-state gas flows. Finally, motivated by the Austrian decarbonization plans for the power sector \cite{EAG2021}, in section~\ref{subsec:Policy} we discuss policy constraints to achieve this and design a green power system constraint with implications for the production of renewable hydrogen.

\subsection{Temporal structure}
\label{subsec:Time}
The LEGO model has a flexible temporal structure, which allows it to represent time either (\textit{i}) by a full chronological time series or (\textit{ii}) by representative periods. Both frameworks allow to specify the temporal resolution, e.g. as hourly, but any resolution can be chosen.
To maintain this level of temporal flexibility, we employ three different temporal indices: chronological periods \textit{p}, representative periods \textit{rp}, and chronological periods within a representative period \textit{k}. In addition, there are two parameters $W^{RP}_{rp}$ and $W^{K}_{k}$ that represent the weights of representative periods $rp$ and the weights (or duration) of periods $k$ respectively. 
First, we explain the full chronological time series case.
Here, all chronological hours $p = 1,2,\dots,8760$ in a year are represented by a single \textit{rp} where each
$k = 1,2,\dots,8760$ is mapped to its corresponding \textit{p}. Since there is only one \textit{rp}, $W^{RP}_{rp}=1$ and since each hour occurs exactly once, $W^{K}_{k}=1$.

In a representative period framework, \textit{rp} and \textit{k} are once again mapped to their corresponding \textit{p}. Mapping and $W^{RP}_{rp}$ are the result of a clustering procedure, e.g. k-medoids.
For the representative periods framework, the following general rules apply: $\sum_{rp} W^{RP}_{rp} = 365$ and $\sum_{rp,k} W^{RP}_{rp} W^{K}_{k} = 8760$.

To give the reader a better understanding, we demonstrate this with a stylized example.
Let us assume one wants to model one year based on five representative days with hourly resolution.
In this case $p = 1,2,\dots,8760$, $rp = 1,2,\dots,5$, and $k = 1,2,\dots,24$.
$W^{K}_{k}=1$, since each hour of the day occurs exactly once.
Furthermore, let us assume that the clustering algorithm determines $W^{RP}_{rp} = 73$ for each \textit{rp}.
This results in $\sum_{rp} W^{RP}_{rp} = 5\cdot73 = 365$ and $\sum_{rp,k} W^{RP}_{rp} W^{K}_{k} = 5\cdot73\cdot1\cdot24 = 8760$.

\subsection{Objective function and bounds}
\label{sec:OF}
The proposed ESM is based in the objective of minimizing total system cost.
In the following, we present the objective function and several generic bounds of the ESM where we focus on the natural gas and hydrogen sectors. In order to give the reader a holistic perspective, we include the essential cost elements of the power sector and summarize the elements of the power sector that are not relevant or original contributions of this paper, e.g. demand-side management, under the term (\textit{xviii}).
For more details on this, the interested reader is referred to \cite{Wogrin2022-LEGO}. All indices, parameters, and variables are described in the Appendix~\ref{sec:Appendix}.

The objective function (\ref{eqn:OF}) includes:
(\textit{i})     cost for supplying natural gas to the system;
(\textit{ii})    operation and maintenance (OM) costs of gas-fired thermal units; 
(\textit{iii})   cost for startup, commitment, and generation of thermal units (except gas-fired thermal units);
(\textit{iv})    OM cost of renewable units; 
(\textit{v})     OM cost of storage units (power system);
(\textit{vi})    cost of electricity non-supplied;
(\textit{vii})   cost of hydrogen and natural gas non-supplied;
(\textit{viii})  CO\textsubscript{2} costs of thermal units (except gas-fired thermal units);
(\textit{ix})    CO\textsubscript{2} costs of gas-fired thermal units;
(\textit{x})     cost of CO\textsubscript{2} emissions in the industry sector;
(\textit{xi})    investment costs for power generation units;
(\textit{xii})   investment costs for transmission lines;
(\textit{xiii})  investment costs for hydrogen units;
(\textit{xiv})   OM costs for hydrogen units;
(\textit{xv})    investment costs for natural gas units;
(\textit{xvi})   OM costs for natural gas units;
(\textit{xvii})  investment costs for hydrogen pipelines; and
(\textit{xviii}) summarized cost elements not relevant or original contributions of this paper.
The system-wide natural gas demand is met from gas wells (\textit{i}). This implies that unit commitment costs (except for OM costs (\textit{ii}) of combined cycle gas turbines (CCGTs) and open cycle gas turbines (OCGTs)), costs for natural gas consumption from steam-methane reforming (SMR) units, and costs of meeting natural gas demand other than for power generation are accounted for implicitly. Furthermore, natural gas and/or hydrogen consumption of compressor units is also accounted for implicitly (see section~\ref{subsec:GasNetwork}).
Constraints (\ref{eqn:PNS}-\ref{eqn:H2/CH4NS}) establish lower and upper bounds for power, hydrogen, and natural gas non-supplied.
Finally, (\ref{eqn:BoundInvProdCons}-\ref{eqn:BoundInvTrans}) establish non-negativity and limit the investments in generation, hydrogen units, natural gas units, and transmission infrastructure (power, natural gas, and hydrogen).
\begin{subequations}
\begin{align}
    min \sum_{rp,k} W^{RP}_{rp} W^{K}_{k} \Big(  \sum_{ch4w}        \underbrace{   C^{CH4}  p^{CH4}_{rp,k,ch4w}}_\textit{i}
                                                +\sum_{t=gas   }    \underbrace{   C^{OM}_t p^{E}_{rp,k,t}}_\textit{ii} \nonumber \\ 
                                                +\sum_{t\neq gas} ( \underbrace{   C^{SU}_t  y_{rp,k,t} 
                                                                                  +C^{UP}_t  u_{rp,k,t} 
                                                                                  +C^{VAR}_t p^{E}_{rp,k,t}}_\textit{iii} )  \nonumber \\
                                                +\sum_r              \underbrace{  C^{OM}_r p^{E}_{rp,k,r}}_\textit{iv}
                                                +\sum_s              \underbrace{  C^{OM}_s p^{E}_{rp,k,s}}_\textit{v} 
                                                +\sum_i              \underbrace{  C^{ENS}  pns_{rp,k,i}}_\textit{vi}  \nonumber \\       
                                                +\sum_{m,cl}        \underbrace{( C^{H2NS} h2ns_{rp,k,m,cl} 
                                                                                  +C^{CH4NS}ch4ns_{rp,k,m,cl})}_\textit{vii} \Big) \nonumber \\
                                                + \sum_{t \neq gas} \underbrace{C^{CO2} E_{t} ( C^{SU}_t y_{rp,k,t} 
                                                                                               +C^{UP}_t u_{rp,k,t}
                                                                                               +C^{VAR}_t p^{E}_{rp,k,t} )}_\textit{viii}                             \nonumber \\
                                                + \sum_{t=gas}\underbrace{C^{CO2} E_{t} (cs^{CH4,E}_{rp,k,t} + cs^{CH4,Aux}_{rp,k,t})}_\textit{ix} \nonumber \\
                                                + \sum_{cls(cl,sec=industry)}\underbrace{C^{CO2} E_{cl} d^{CH4}_{rp,k,m,cl}}_\textit{x}
                                                \Big) \nonumber\\
   + \sum_g \underbrace{C^{INV}_g x_g}_\textit{xi}
   + \sum_{ijcc(i,j,c)} \underbrace{C^{L,Inv}_{i,j,c} x^L_{i,j,c}}_\textit{xii} \nonumber\\
   + \sum_{h2u} \underbrace{(C^{INV}_{h2u} x^{H2}_{h2u})}_\textit{xiii}
              + \underbrace{(C^{OM}_{h2u} (x^{H2}_{h2u}+EU^{H2}_{h2u}))}_\textit{xiv} \nonumber \\
   + \sum_{ch4u} \underbrace{(C^{INV}_{ch4u} x^{CH4}_{ch4u})}_\textit{xv}
               + \underbrace{(C^{OM}_{ch4u} (x^{CH4}_{ch4u}+EU^{CH4}_{ch4u}))}_\textit{xvi} 
               \nonumber \\
   + \sum_{mnlc(m,n,l)} \underbrace{C^{Pipe,Inv}_{m,n,l} x^{Pipe}_{m,n,l}}_\textit{xvii}  
   + \underbrace{other}_\textit{xviii} \label{eqn:OF}\\
   0 \leq pns_{rp,k,i} \leq D^E_{rp,k,i} \quad \forall rp,k,i \label{eqn:PNS} \\
   0 \leq h2ns_{rp,k,m,cl}, ch4ns_{rp,k,m,cl} \leq M \quad \forall rp,k,m,cl \label{eqn:H2/CH4NS} \\
   x_{g}, x^{H2}_{h2u}, x^{CH4}_{ch4u} \in \mathbb{Z}^{+,0}, \nonumber\\
            x_{g}        \leq \overline{X}_{g},
            x^{H2}_{h2u} \leq \overline{X}^{H2}_{h2u},
            x^{CH4}_{ch4u} \leq \overline{X}^{CH4}_{ch4u}
            \quad \forall g, h2u, ch4u \label{eqn:BoundInvProdCons} \\
   x^L_{i,j,c}, x^{Pipe}_{m,n,l} \in \{0,1\},
            x^L_{i,j,c}      \leq \overline{X}^{L}_{i,j,c},
            x^{Pipe}_{m,n,l} \leq \overline{X}^{Pipe}_{m,n,l} \nonumber \\ \quad \forall ijcc(i,j,c), mnlc(m,n,l) \label{eqn:BoundInvTrans} 
\end{align}
\end{subequations}

\subsection{Natural gas sector}
\label{subsec:NG-sector}
This section introduces the mathematical framework governing the supply and demand side of the natural gas sector.
In the presented ESM the entire natural gas demand is supplied by natural gas wells (\ref{eqn:Gas-Wells}).
\begin{align}
  0 \leq p^{CH4}_{rp,k,ch4w} \leq \overline{P}_{ch4w} (x^{CH4}_{ch4w}+EU^{CH4}_{ch4w}) \quad \forall rp,k,ch4w \label{eqn:Gas-Wells}
\end{align}
The gas demand side is governed by (\ref{eqn:Gas-DemBlend}).
The framework allows to partially substitute natural gas demand in various sectors, e.g. iron and steel, chemistry, households etc. by blending hydrogen.
To this end, constraint (\ref{eqn:Gas-Blend-Demand}) ensures that the energy content of the resulting blend is sufficient, where natural gas and hydrogen demand (multiplied by their lower heating values) are additive. This gives an accuracy of +0.5\% for the blend's heating value \cite{Cerbe2016-2}. Furthermore, we assume that the resulting hydrogen demand has to follow the same temporal pattern as the original natural gas demand.
Constraints (\ref{eqn:Gas-BoundsBlend-CH4Demand}-\ref{eqn:Gas-BoundsBlend-H2Demand}) establish lower and upper bounds for the natural gas variable and the volumetric substitution of hydrogen per sector respectively.
\begin{subequations}
\label{eqn:Gas-DemBlend}
\begin{align}
  D^{Gas}_{rp,k,m,cl} H^{CH4} = d^{CH4}_{rp,k,m,cl} H^{CH4} + d^{H2}_{rp,k,m,cl} H^{H2} 
  \quad \forall rp,k,m,cl \label{eqn:Gas-Blend-Demand} \\
  0 \leq d^{CH4}_{rp,k,m,cl} \leq D^{Gas}_{rp,k,m,cl} \quad \forall rp,k,m,cl \label{eqn:Gas-BoundsBlend-CH4Demand} \\
  \underline{SR}_{cl}^{H2} d^{CH4}_{rp,k,m,cl} \leq d^{H2}_{rp,k,m,cl} \leq \overline{SR}_{cl}^{H2} d^{CH4}_{rp,k,m,cl} \quad \forall rp,k,m,cl \label{eqn:Gas-BoundsBlend-H2Demand}
\end{align}
\end{subequations}

Constraints (\ref{eqn:Gas-thermals}) represent gas-fired CCGTs and OCGTs. In reality, the relation between gas consumption and power output of a thermal generation unit, referred to a heat rate, is non-linear. However, in the literature commonly a linear relation is assumed, e.g. \cite{Jia2021,Zhang2015a}.
Equation (\ref{eqn:Gas-Conversion}) establishes the relation for the conversion of gas to power.
In addition, equation (\ref{eqn:Gas-Auxiliary}) includes the auxiliary gas consumption associated with the startup and commitment of the unit.
We formulate the total gas consumption in two separate equations (\ref{eqn:Gas-Conversion}-\ref{eqn:Gas-Auxiliary}) to preserve information about the resource (natural gas or hydrogen) used for power generation (see section~\ref{subsec:Policy}).
The remaining constraints are:
lower and upper bounds for natural gas consumption (\ref{eqn:Gas-ConsPower-GasBounds}-\ref{eqn:Gas-ConsAux-GasBounds});
lower and upper bound for hydrogen blending (\ref{eqn:Gas-BoundsBlendPower-H2GT}-\ref{eqn:Gas-BoundsBlendAux-H2GT}); and
lower and upper bound for power output (\ref{eqn:Gas-BoundPowerOutput-GT}).
\begin{subequations}
\label{eqn:Gas-thermals}
\begin{align}
  cs^{CH4,E}_{rp,k,t} H^{CH4} + cs^{H2,E}_{rp,k,t} H^{H2} = p^{E}_{rp,k,t} CS^{V}_t  \quad \forall rp,k,t=gas \label{eqn:Gas-Conversion} \\
  cs^{CH4,Aux}_{rp,k,t} H^{CH4} + cs^{H2,Aux}_{rp,k,t} H^{H2} = y_{rp,k,t} CS^{SU}_t / W^{K}_{k} + u_{rp,k,t} CS^{UP}_t \nonumber \\ \quad \forall rp,k,t=gas \label{eqn:Gas-Auxiliary} \\
  0 \leq cs^{CH4,E}_{rp,k,t} \leq CS^{V}_t \overline{P}^{E}_{t} / H^{CH4} (x_t + EU_t) \quad \forall rp,k,t=gas \label{eqn:Gas-ConsPower-GasBounds} \\
  0 \leq cs^{CH4,Aux}_{rp,k,t} \leq (y_{rp,k,t} CS^{SU}_t / W^{K}_{k} + u_{rp,k,t} CS^{UP}_t) / H^{CH4} \quad \forall rp,k,t=gas \label{eqn:Gas-ConsAux-GasBounds} \\
  \underline{B}^{H2} cs^{CH4,E}_{rp,k,t} \leq cs^{H2,E}_{rp,k,t} \leq \overline{B}^{H2} cs^{CH4,E}_{rp,k,t} \quad \forall rp,k,t=gas \label{eqn:Gas-BoundsBlendPower-H2GT} \\
  \underline{B}^{H2} cs^{CH4,Aux}_{rp,k,t} \leq cs^{H2,Aux}_{rp,k,t} \leq \overline{B}^{H2} cs^{CH4,Aux}_{rp,k,t} \quad \forall rp,k,t=gas \label{eqn:Gas-BoundsBlendAux-H2GT} \\
  0 \leq p^{E}_{rp,k,t} \leq \overline{P}^{E}_t (x_t + EU_t) \quad \forall rp,k,t=gas \label{eqn:Gas-BoundPowerOutput-GT}
\end{align}
\end{subequations}

Finally, the model also includes long-term natural gas storage units. Their mathematical formulation is basically identical to long-term hydrogen storage units presented in (\ref{eqn:H2-Storage}).

\subsection{Hydrogen sector}
\label{subsec:H2-sector}
In section~\ref{subsec:H2-sector}, we introduce the utilities of the (future) hydrogen sector. 
This includes detailed formulations for electrolyzer (EL), steam-methane reforming (SMR),(long-term) hydrogen storage, and fuel cell (FC) units.

EL units consume electricity to produce hydrogen by the process of water electrolysis (\ref{eqn:H2-Conversion-EL}).
Constraints (\ref{eqn:H2-Bounds-Cons-EL}-\ref{eqn:H2-Bounds-Prod-EL}) establish lower and upper bounds for electricity consumption and hydrogen production respectively.
\begin{subequations}
\label{eqn:H2-EL}
\begin{align}
   p^{H2}_{rp,k,h2g} = cs^{E}_{rp,k,h2g} HPE_{h2g} \quad \forall rp,k,h2g \label{eqn:H2-Conversion-EL} \\
  0 \leq cs^{E}_{rp,k,h2g} \leq \overline{P}^E_{h2g} (x^{H2}_{h2g} + EU^{H2}_{h2g}) \quad \forall rp,k,h2g \label{eqn:H2-Bounds-Cons-EL} \\
  0 \leq p^{H2}_{rp,k,h2g} \leq \overline{P}^E_{h2g} HPE_{h2g} (x^{H2}_{h2g} + EU^{H2}_{h2g}) \quad \forall rp,k,h2g \label{eqn:H2-Bounds-Prod-EL}
\end{align}
\end{subequations}

Today, around 75\% of the global hydrogen demand is supplied from the process of steam-methane reforming \cite{Gul2019}. Unlike water electrolysis, SMR produces hydrogen from natural gas (\ref{eqn:H2-Conversion-SMR}), which serves as a fuel for steam generation and as feedstock for the process. In the formulation, this is accounted for in the hydrogen per natural gas ratio (or efficiency factor). Similarly to EL units, constraints (\ref{eqn:H2-Bounds-Cons-SMR}-\ref{eqn:H2-Bounds-Prod-SMR}) establish lower and upper bounds for natural gas consumption and hydrogen production.
\begin{subequations}
\begin{align}
  p^{H2}_{rp,k,h2p} = cs^{CH4}_{rp,k,h2p} HPC_{h2p} \quad \forall rp,k,h2p \label{eqn:H2-Conversion-SMR} \\
  0 \leq cs^{CH4}_{rp,k,h2p} \leq (\overline{P}^{H2}_{h2p}/HPC_{h2p}) (x^{H2}_{h2p} + EU^{H2}_{h2p}) \quad \forall rp,k,h2p \label{eqn:H2-Bounds-Cons-SMR} \\
  0 \leq p^{H2}_{rp,k,h2p} \leq \overline{P}^{H2}_{h2p} (x^{H2}_{h2p} + EU^{H2}_{h2p}) \quad \forall rp,k,h2p \label{eqn:H2-Bounds-Prod-SMR}
\end{align}
\end{subequations} 

Constraints (\ref{eqn:H2-SIntraSOCDef}-\ref{eqn:H2-Bounds-Ch}) describe hydrogen storage units based on their state of charge (SOC). Depending on the temporal framework one chooses (see section~\ref{subsec:Time}), modeling hydrogen storage units relies on a different concept, which we describe in the following lines.

In models that represent time chronologically, e.g. all hours of a year, both short- and long-term storage units can be modeled based on a common set of SOC constraints. 
In a representative period framework (which is common for infrastructure planning due to the implied computational complexity), a set of common SOC constraints gives an accurate representation of (short-term) storage units within the representative period, e.g. a day. However, accurate representation of long-term storage units in a representative periods framework is difficult.
In order to model both short- and long-term hydrogen and natural gas storage units, we separate the SOC constraints into inter-period and intra-period constraints.

Short-term hydrogen storage units, e.g. steel tanks, are governed by an intra-period SOC concept (\ref{eqn:H2-SIntraSOCDef}-\ref{eqn:H2-SUBBound-IntraSOC}). In this concept, the last hour and the first hour of each representative period are linked together, which establishes a cyclic relation. In other words, $k24$ of representative day $rp1$ is followed by $k1$ of representative day $rp1$. This cyclic relation prevents depletion of short-term hydrogen storage units towards the end of an representative period and is indicated by the double minus notation, e.g. $k\!-\!-1$.

The typical characteristic (and also the purpose) of a short-term storage unit is its ability to cycle (completely deplete and recharge the stored energy) frequently, eventually even several times a day.
In contrast, long-term hydrogen storage units, e.g. salt caverns, are designed to serve a different purpose, e.g. inter-seasonal hydrogen storage. To capture their operational decisions and associated long-term effects in a representative period framework, long-term hydrogen storage units follow a inter-period SOC concept (\ref{eqn:H2-SInterSOCDef}-\ref{eqn:H2-SUBInterSOC}) \cite{Tejada-Arango2018}.
This concept is based on a moving window (MOW) that moves along the index \textit{p} and imposes the inter-period SOC constraints at each multiple of the specified duration of the MOW. In the process, production and consumption decisions that occur within the MOW are accounted for.

In a full chronological time series framework, which corresponds to considering a single representative period, e.g. of 8760 hours, all storage units are represented only by the intra-period SOC formulation and the cyclic relation described above is omitted. Constraints (\ref{eqn:H2-Bounds-Dis}-\ref{eqn:H2-Bounds-Ch}) apply for all of the described cases and establish lower and upper bounds on the production and consumption of hydrogen storage units respectively.
\begin{subequations}
\label{eqn:H2-Storage}
\begin{align}
  intra^{H2}_{rp,k,h2s} = intra^{H2}_{rp,k--1,h2s} - p^{H2}_{rp,k,h2s} W^K_k / \eta^{DIS}_{h2s} + cs^{H2}_{rp,k,h2s} W^K_k \eta^{CH}_{h2s} \nonumber \\ \quad \forall rp,k,h2s \label{eqn:H2-SIntraSOCDef} \\
 intra^{H2}_{rp,k,h2s} \geq \underline{R}^{H2}_{h2s} \overline{P}^{H2}_{h2s} ETP_{h2s} (x^{H2}_{h2s} + EU^{H2}_{h2s}) \quad \forall rp,k,h2s  \label{eqn:H2-SLBBound-IntraSOC} \\
 intra^{H2}_{rp,k,h2s} \leq \overline{P}^{H2}_{h2s} ETP_{h2s} (x^{H2}_{h2s} + EU^{H2}_{h2s})  \quad \forall rp,k,h2s  \label{eqn:H2-SUBBound-IntraSOC}\\
  inter^{H2}_{p,h2s} = inter^{H2}_{p-MOW,h2s} + InRes^{H2}_{h2s,p=MOW} (x^{H2}_{h2s} + EU^{H2}_{h2s})  \nonumber \\
  + \sum_{\Gamma (p-MOW \leq pp\leq p,rp,k)} (- p^{H2}_{ rp,k,h2s} W_k^K / \eta^{DIS}_{h2s} 
                                            + cs^{H2}_{rp,k,h2s} W_k^K   \eta^{CH}_{h2s})
                                             \quad \forall p,h2s  
                                             \label{eqn:H2-SInterSOCDef} \\
  inter^{H2}_{p,h2s} \geq \underline{R}^{H2}_{h2s} \overline{P}^{H2}_{h2s} ETP_{h2s} (x^{H2}_{h2s} + EU^{H2}_{h2s})  \quad \forall h2s, p: mod(p,MOW)=0 \label{eqn:H2-SLBInterSOC} \\
  inter^{H2}_{p,h2s} \leq \overline{P}^{H2}_{h2s} ETP_{h2s} (x^{H2}_{h2s} + EU^{H2}_{h2s})  \quad \forall h2s, p: mod(p,MOW)=0  \label{eqn:H2-SUBInterSOC} \\
  inter^{H2}_{p,h2s} = InRes^{H2}_{h2s} (x^{H2}_{h2s} + EU^{H2}_{h2s}) \quad \forall h2s, p=CARD(p)  \label{eqn:H2-SCyclicInter} \\
  0 \leq p^{H2}_{rp,k,h2s} \leq \overline{P}^{H2}_{h2s} (x^{H2}_{h2s} + EU^{H2}_{h2s}) \quad \forall rp,k,h2s  \label{eqn:H2-Bounds-Dis}\\
  0 \leq cs^{H2}_{rp,k,h2s} \leq \overline{CS}^{H2}_{h2s} (x^{H2}_{h2s} + EU^{H2}_{h2s}) \quad \forall rp,k,h2s  \label{eqn:H2-Bounds-Ch}
\end{align}
\end{subequations}

Finally, FC units reverse the electrolysis process and generate power from hydrogen (\ref{eqn:H2-Conversion-FC}).
Constraints (\ref{eqn:H2-Bounds-Cons-FC}-\ref{eqn:H2-Bounds-Prod-FC}) establish lower and upper bounds for hydrogen consumption and power generation respectively.
\begin{subequations}
\label{eqn:H2-FC}
\begin{align}
   p^{E}_{rp,k,h2f} = cs^{H2}_{rp,k,h2f} EPH_{h2f} \quad \forall rp,k,h2f \label{eqn:H2-Conversion-FC} \\
  0 \leq cs^{H2}_{rp,k,h2f} \leq \overline{P}^{H2}_{h2f} (x^{H2}_{h2f} + EU^{H2}_{h2f}) \quad \forall rp,k,h2f \label{eqn:H2-Bounds-Cons-FC} \\
  0 \leq p^{E}_{rp,k,h2f} \leq \overline{P}^E_{h2f} EPH_{h2f} (x^{H2}_{h2f} + EU^{H2}_{h2f}) \quad \forall rp,k,h2f \label{eqn:H2-Bounds-Prod-FC}
\end{align}
\end{subequations}

\subsection{Gas network model}
\label{subsec:GasNetwork}
The gas transmission network connects gas production facilities, e.g. natural gas wells, and gas demand. 
Today, it is already possible to transport natural gas-hydrogen blends via the gas transmission network - at least for small blending rates, e.g. 10\% in Austria \cite{GasConnect2021}.
In the future it is envisaged to establish a dedicated European hydrogen network \cite{Jens2021}. In reality, gas transmission systems comprise pipelines, measuring equipment, equipment for gas cleaning, compressor units (CUs), cooling sections and gas drying units. In ESMs, however, generally only pipelines and - depending on the level of detail and the gas flow formulation - CUs are considered.

\begin{figure}[h]
\centerline{\includegraphics[scale=0.40]{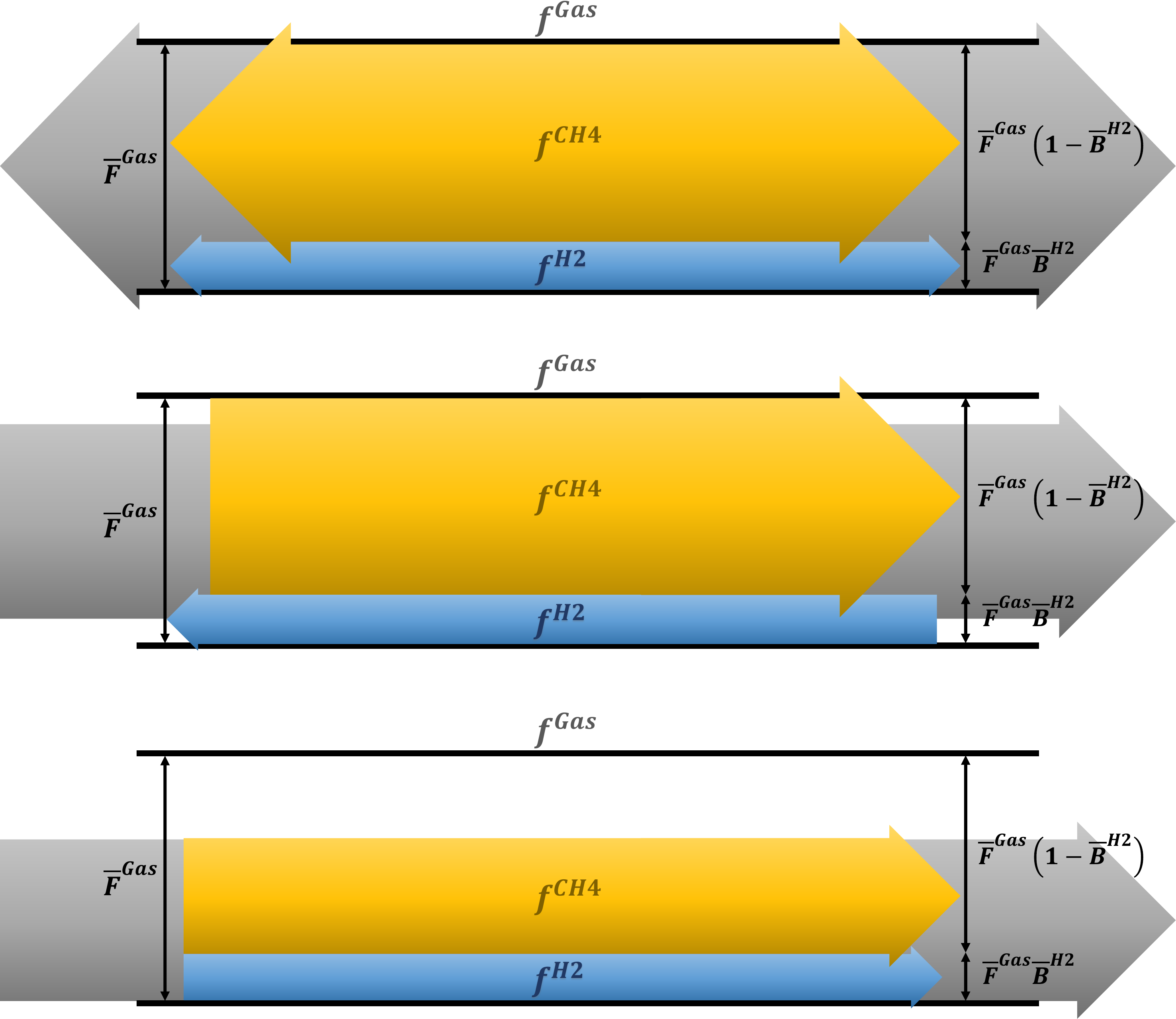}}
\caption{Graphical representation of pipeline gas flow based on the S-TP formulation in the context of natural gas and hydrogen blending. Top: Flow directions can be altered in each time step (\textit{i}); Center: Natural gas and hydrogen flow in opposite directions (\textit{ii}); Bottom: Violation of maximum volumetric blending rate (\textit{iii}).}
\label{fig:S-TP}
\end{figure}

For gas network modeling the formulation of gas flow is the centerpiece \cite{Klatzer2022}.
As outlined in section~\ref{sec:Literature}, gas flow formulations can comprise various levels of technical detail.
The most basic approach is to disregard the physical laws governing gas flows and consider a standard transport problem (S-TP), which is a linear program. Despite its simplicity the S-TP is a common approach in the literature, e.g. \cite{Nunes2018,Zhao2018}.
However, re-formulating and applying the S-TP in the context of natural gas and hydrogen blending (\ref{eqn:S-TP_Gas_Flow}) has significant implications for investment decisions in the power, natural gas, and hydrogen sectors, which can lead to sub-optimal planning and ultimately non-supplied hydrogen \cite{KlatzerGM2023}.
The reasons for this are: (\textit{i}) With the S-TP, flow directions can be altered in each time step, (\textit{ii}) it allows natural gas and hydrogen to flow in opposite directions in a pipeline, which is not possible in reality, and (\textit{iii}) it cannot guarantee compliance with the maximum hydrogen blending rate, e.g. 10\% of the actual natural gas flow. Fig.~\ref{fig:S-TP} provides a graphical representation of the issues stated.
\begin{subequations}
\label{eqn:S-TP_Gas_Flow} 
\begin{align}
  f^{Gas}_{rp,k,m,n,l} = f^{CH4}_{rp,k,m,n,l} + f^{H2}_{rp,k,m,n,l}  \quad \forall rp,k,mnl(m,n,l) \label{eqn:S-TP-1} \\
 -\overline{F}^{Gas}_{m,n,l} \overline{B}^{H2} \leq f^{H2}_{rp,k,m,n,l} \leq \overline{F}^{Gas}_{m,n,l} \overline{B}^{H2} \quad  \forall rp,k,mnl(m,n,l) \label{eqn:S-TP-2} \\
  -\overline{F}^{Gas}_{m,n,l} (1-\overline{B}^{H2}) \leq f^{CH4}_{rp,k,m,n,l} \leq \overline{F}^{Gas}_{m,n,l} (1-\overline{B}^{H2}) \quad  \forall rp,k,mnl(m,n,l) \label{eqn:S-TP-3} \\
 -\overline{F}^{Gas}_{m,n,l} x^{Pipe}_{m,n,l} \leq f^{Gas}_{rp,k,m,n,l} \leq  
  \overline{F}^{Gas}_{m,n,l} x^{Pipe}_{m,n,l} 
  \quad \forall rp,k,mnlc(m,n,l) \label{eqn:S-TP-4}
\end{align}
\end{subequations}

\begin{figure}[h]
\centerline{\includegraphics[scale=0.40]{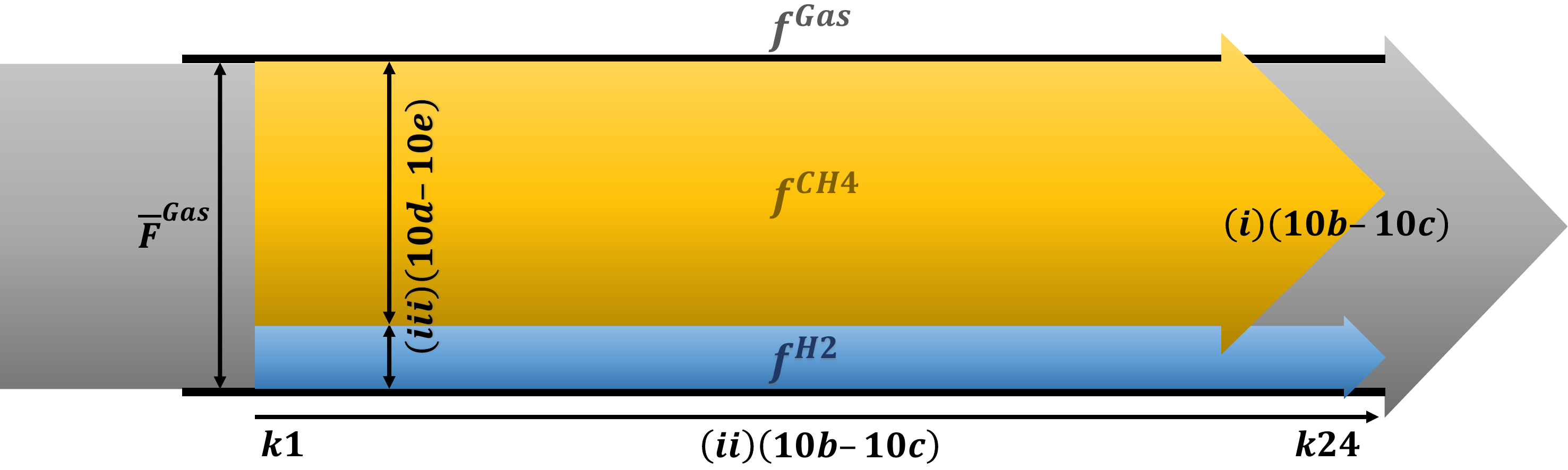}}
\caption{Graphical representation of pipeline gas flow based on the B-TP formulation. Constraints (\ref{eqn:H2FlowDir-TP}-\ref{eqn:CH4FlowDir-TP}) ensure that (\textit{i}) natural gas and hydrogen flow in the same direction in a pipeline and (\textit{ii}) that gas flows are subject to one decision per representative day. Constraints (\ref{eqn:Gas-ExiLB-H2-TP}-\ref{eqn:Gas-ExiUB-H2-TP}) ensure that (\textit{iii}) the volumetric blending rate is not violated.}
\label{fig:B-TP}
\end{figure}

To overcome the issues of the S-TP, we proposed a novel blending transport problem (B-TP) \cite{KlatzerGM2023} (\ref{eqn:B-TP_Gas_Flow}). In contrast to the S-TP, the B-TP is based on a MILP framework. Therein, the binary variable $\alpha^{Gas}_{rp,m,n,l}$ ensures that (\textit{i}) natural gas and hydrogen flow in the same direction in a pipeline, (\textit{ii}) that gas flows are subject to one decision per representative day, which is a good approximation of reality, and (\textit{iii}) that the volumetric blending rate is not violated. Fig.~\ref{fig:B-TP} provides a graphical representation of the B-TP formulation.
\begin{subequations}
\label{eqn:B-TP_Gas_Flow} 
\begin{align}
 f^{Gas}_{rp,k,m,n,l} = f^{CH4}_{rp,k,m,n,l} + f^{H2}_{rp,k,m,n,l}  \quad \forall rp,k,mnl(m,n,l) \label{eqn:Def-GasFlowVar-TP} \\
 (\alpha^{Gas}_{rp,m,n,l} - 1) M \leq f^{H2}_{rp,k,m,n,l} \leq 
  \alpha^{Gas}_{rp,m,n,l} M
  \quad \forall rp,k,mnl(m,n,l) \label{eqn:H2FlowDir-TP}\\ 
 (\alpha^{Gas}_{rp,m,n,l} - 1) M \leq f^{CH4}_{rp,k,m,n,l} \leq
  \alpha^{Gas}_{rp,m,n,l} M
  \quad \forall rp,k,mnl(m,n,l) \label{eqn:CH4FlowDir-TP}\\
  f^{H2}_{rp,k,m,n,l} \geq -\alpha^{Gas}_{rp,m,n,l} M + \overline{B}^{H2} f^{CH4}_{rp,k,m,n,l}
  \quad \forall rp,k,mnl(m,n,l) \label{eqn:Gas-ExiLB-H2-TP} \\
 (1-\alpha^{Gas}_{rp,m,n,l}) M + \overline{B}^{H2} f^{CH4}_{rp,k,m,n,l} \geq f^{H2}_{rp,k,m,n,l}
 \quad \forall rp,k,mnl(m,n,l) \label{eqn:Gas-ExiUB-H2-TP} \\
 -\overline{F}^{Gas}_{m,n,l} \leq f^{Gas}_{rp,k,m,n,l} \leq \overline{F}^{Gas}_{m,n,l} \quad  \forall rp,k,mnl(m,n,l) \label{eqn:Gas-ExiBounds-Gas-TP} \\
 -\overline{F}^{Gas}_{m,n,l} x^{Pipe}_{m,n,l} \leq f^{Gas}_{rp,k,m,n,l} \leq  
  \overline{F}^{Gas}_{m,n,l} x^{Pipe}_{m,n,l} 
  \quad \forall rp,k,mnlc(m,n,l) \label{eqn:Gas-CanBounds-Gas-TP}
\end{align}
\end{subequations}

In this paper, we extend our approach in order to capture the physical relationship between gas flow and gas pressure in high-pressure pipelines, which can be described by the steady-state gas flow equation \cite{Mischner2020}.
This relationship, however, is non-linear and non-convex, which is problematic in the context of (large-scale) energy system modeling. To this end, a variety of linearization methods have been applied. This includes Taylor series \cite{Ordoudis2019}, Newton-Raphson method, e.g. used for commercial software tools like PSS SINCAL and NEPLAN \cite{Kriechbaum2018}, iterative methods \cite{Lohr2020}, and piecewise linearization \cite{Correa-Posada2015,Zhang2018,Jia2021}.

\begin{figure}[h]
\centerline{\includegraphics[scale=0.59]{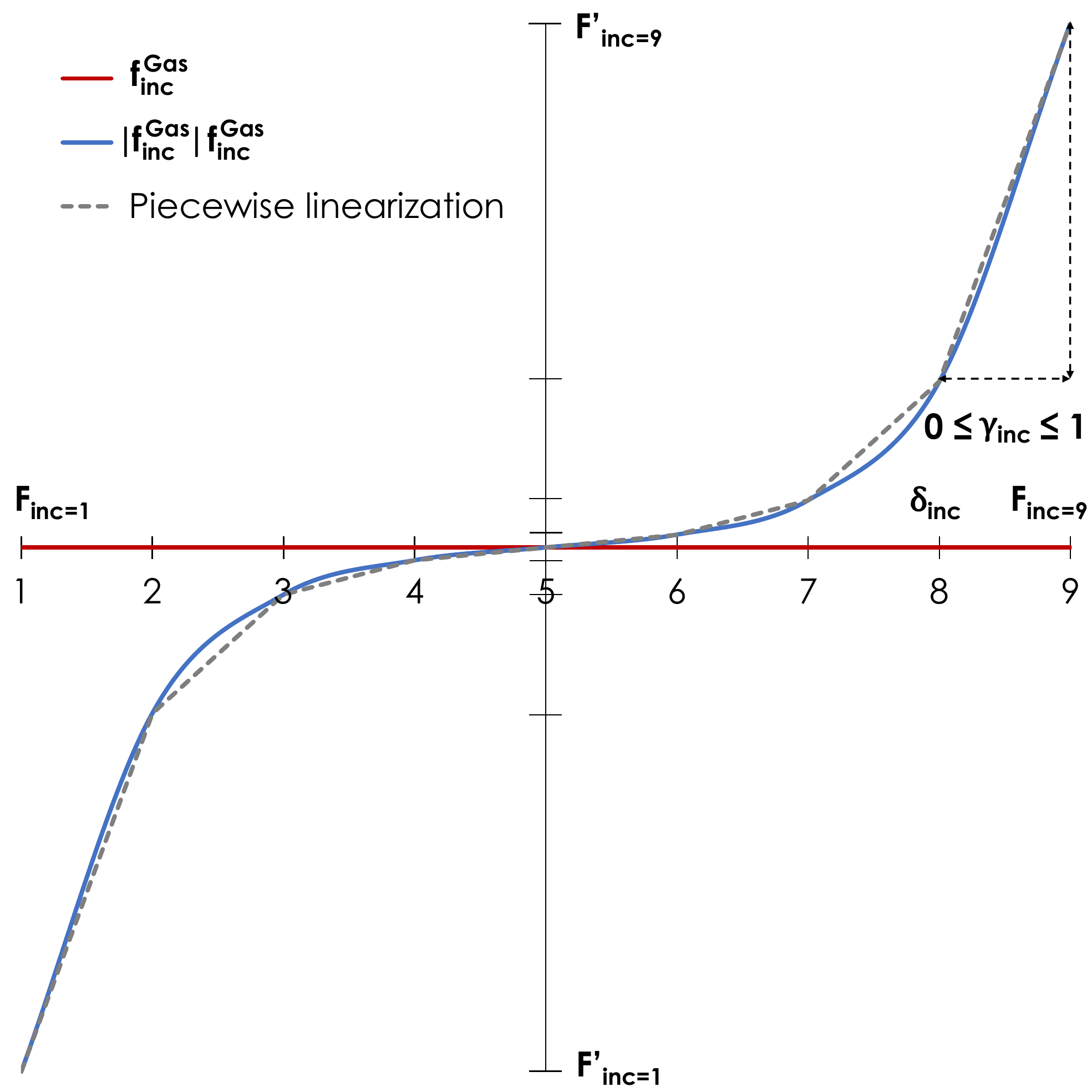}}
\caption{Graphical representation of the incremental piecewise linearization method.}
\label{fig:B-PP}
\end{figure}

The novel formulation to model gas flows based on the steady-state gas flow equation proposed in this paper extends the work of \cite{Correa-Posada2015} to include natural gas and hydrogen blending.
For linearization, we adopt an incremental (INC) piecewise linearization method (MILP framework), which has been found to computationally outperform other piecewise linearization methods \cite{Correa-Posada2014}.
For a better understanding we want to briefly explain the underlying concept of the linearization method presented in (\ref{eqn:GasFlow-INC}-\ref{eqn:GasFlow-Bounds-Press-INC}).
Equation (\ref{eqn:GasFlow-INC}) describes the relationship between the directional quadratic gas flow $|f^{Gas}_{rp,k,m,n,l}|f^{Gas}_{rp,k,m,n,l}$ expressed as a linear combination of function values $F'_{inc,m,n,l}$ (LHS) and the squared pressure variables at the start and end-point of a pipeline (RHS).
The continuous variable $\gamma_{rp,k,inc,m,n,l}$ links (\ref{eqn:GasFlow-INC}) and (\ref{eqn:Flow-Variable-INC}) and thus establishes the relationship between the linear gas flow $f^{Gas}_{rp,k,m,n,l}$ and the squared pressure variables.
Finally, constraint (\ref{eqn:Filling-Cond-INC}) utilizes the binary variable $\delta_{rp,k,inc,m,n,l}$ to ensure that the order of gas flow increments is preserved.
Fig.~\ref{fig:B-PP} provides a graphical representation of the incremental piecewise linearization method.

On that basis, we extend the INC formulation to include pipeline transmission expansion planning (TEP) and blending of natural gas and hydrogen.
For TEP, we introduce the slack variable $\rho_{rp,k,m,n,l}$ in (\ref{eqn:Flow-Variable-INC}) and add constraints (\ref{eqn:Bounds-ExiGasSlack-INC}-\ref{eqn:GasFlow-TEP-INC}). 
The slack variable ensures that the equality in (\ref{eqn:Flow-Variable-INC}) can be met in case a specified candidate pipeline is not built. 
To model natural gas and hydrogen blending, we introduce constraints (\ref{eqn:Def-GasFlowVar-INC}-\ref{eqn:Gas-Bounds-Dir-INC}). As with the B-TP, the binary variable $\alpha^{Gas}_{rp,m,n,l}$ ensures coherence of natural gas and hydrogen flow direction (\ref{eqn:H2FlowDir-INC}-\ref{eqn:CH4FlowDir-INC}) and limits the number of decisions on gas flow direction to one per representative period.
Finally, constraints (\ref{eqn:Gas-ExiLB-H2-INC}-\ref{eqn:Gas-ExiUB-H2-INC}) establish lower and upper bounds on the blending rate, which is expressed as a percentage of the actual natural gas flow.
\begin{subequations}
\label{eqn:Inc_Gas_Flow}
\begin{align}
  F'_{inc=1,m,n,l} +\sum^{INC-1}_{inc} \left(F'_{inc+1,m,n,l} - F'_{inc,m,n,l}\right)
  \gamma_{rp,k,inc,m,n,l}  \nonumber\\
  = R^{Gas}_{m,n,l} (p^{sqr}_{rp,k,m} - p^{sqr}_{rp,k,n}) \quad \forall rp,k,mnl(m,n,l) 
  \label{eqn:GasFlow-INC}\\
  F_{inc=1,m,n,l} +\sum^{INC-1}_{inc} \left(F_{inc+1,m,n,l} - F_{inc,m,n,l}\right) 
  \gamma_{rp,k,inc,m,n,l} \nonumber\\
  = \rho_{rp,k,m,n,l} + f^{Gas}_{rp,k,m,n,l} \quad \forall rp,k,mnl(m,n,l)
  \label{eqn:Flow-Variable-INC} \\
\gamma_{rp,k,inc+1,m,n,l} \leq \delta_{rp,k,inc,m,n,l} \leq \gamma_{rp,k,inc,m,n,l} \nonumber\\
\quad \forall rp,k,mnl(m,n,l),inc \in INC\!-\!1 \label{eqn:Filling-Cond-INC} \\
  0 \leq \gamma_{rp,k,inc,m,n,l} \leq 1 \quad \forall rp,k,inc,mnl(m,n,l) \label{eqn:Def-FillVar-INC}\\
  \delta_{rp,k,inc,m,n,l} \in \{0,1\} \quad \forall rp,k,inc,mnl(m,n,l) \label{eqn:Def-FillVar-INC} \\
  0 \leq p^{sqr}_{rp,k,m}, p^{sqr}_{rp,k,n} \leq \overline{P}_{m}^{sqr}, \overline{P}_{n}^{sqr} \quad \forall rp,k,m,n \label{eqn:GasFlow-Bounds-Press-INC} \\
 \rho_{rp,k,m,n,l} = 0 \quad \forall rp,k,mnle(m,n,l) \label{eqn:Bounds-ExiGasSlack-INC}\\
 -(1-x^{Pipe}_{m,n,l}) \overline{F}^{Gas}_{m,n,l} \leq \rho_{rp,k,m,n,l} \leq 
  (1-x^{Pipe}_{m,n,l}) \overline{F}^{Gas}_{m,n,l} \nonumber\\
  \quad \forall rp,k,mnlc(m,n,l) \label{eqn:Bounds-CanGasSlack-INC}\\
 -\overline{F}^{Gas}_{m,n,l} x^{Pipe}_{m,n,l} \leq f^{Gas}_{rp,k,m,n,l} \leq 
  \overline{F}^{Gas}_{m,n,l} x^{Pipe}_{m,n,l} \quad \forall rp,k,mnlc(m,n,l) \label{eqn:GasFlow-TEP-INC}\\
  f^{Gas}_{rp,k,m,n,l} = f^{CH4}_{rp,k,m,n,l} + f^{H2}_{rp,k,m,n,l}
  \quad \forall rp,k,mnl(m,n,l) \label{eqn:Def-GasFlowVar-INC}\\
 (\alpha^{Gas}_{rp,m,n,l} - 1) M \leq f^{H2}_{rp,k,m,n,l} \leq 
  \alpha^{Gas}_{rp,m,n,l} M
  \quad \forall rp,k,mnl(m,n,l) \label{eqn:H2FlowDir-INC}\\ 
 (\alpha^{Gas}_{rp,m,n,l} - 1) M \leq f^{CH4}_{rp,k,m,n,l} \leq
  \alpha^{Gas}_{rp,m,n,l} M
  \quad \forall rp,k,mnl(m,n,l) \label{eqn:CH4FlowDir-INC}\\
  f^{H2}_{rp,k,m,n,l} \geq -\alpha^{Gas}_{rp,m,n,l} M + \overline{B}^{H2} f^{CH4}_{rp,k,m,n,l}
  \quad \forall rp,k,mnl(m,n,l) \label{eqn:Gas-ExiLB-H2-INC} \\
 (1 -\alpha^{Gas}_{rp,m,n,l}) M + \overline{B}^{H2} f^{CH4}_{rp,k,m,n,l} \geq f^{H2}_{rp,k,m,n,l}
 \quad \forall rp,k,mnl(m,n,l) \label{eqn:Gas-ExiUB-H2-INC} \\
\alpha^{Gas}_{rp,m,n,l} \in \{0,1\} \quad \forall rp, mnl(m,n,l) \label{eqn:Gas-Bounds-Dir-INC}
\end{align}
\end{subequations}

Besides pipelines, compressor units (CUs) play an important role in the gas system. They are typically sited in 100-200~km intervals to compensate the pressure drop caused by friction of gas molecules.
When modelling (large-scale) energy systems, it is common to represent the relationship between inlet and outlet pressure of a CU (\ref{eqn:Comp-Ratio-Rel}) as a linear relationship \cite{Klatzer2022}. In addition, the absolute pressure increase of a CU (\ref{eqn:Comp-Ratio-Abs}) is also limited.
Constraint (\ref{eqn:Gas-BoundsCmp-TP}) establishes lower and upper bounds on the combined natural gas and hydrogen flow through a CU. Note that compressor flows are defined as positive variables. 
This ensures that the gas consumption for the compression work of CUs is positive (see (\ref{eqn:Balances-INC})).
Finally, (\ref{eqn:H2-BoundsCmp-TP}) establishes lower and upper bounds on the hydrogen flow through a CU. Since hydrogen has a lower density than natural gas, this is crucial for the sound operation of CUs \cite{Haeseldonckx2007}.
\begin{subequations}
\label{eqn:Compressors}
\begin{align}
  p^{sqr}_{rp,k,n} \leq \Lambda^{sqr}_{m,n,l} p^{sqr}_{rp,k,m} \quad \forall rp,k,cmp(m,n,l) \label{eqn:Comp-Ratio-Rel}\\
  0 \leq p^{sqr}_{rp,k,n} - p^{sqr}_{rp,k,m} \leq \overline{P}_{m}^{sqr} - \bigg(\sqrt{\overline{P}_{m}^{sqr}} - \overline{\Lambda}_{m,n,l}\bigg)^2 \quad \forall rp,k,cmp(m,n,l) \label{eqn:Comp-Ratio-Abs} \\
  0 \leq f^{Cmp,CH4}_{rp,k,m,n,l} + f^{Cmp,H2}_{rp,k,m,n,l} \leq \overline{F}^{Cmp,Gas}_{m,n,l}
  \quad \forall rp,k,cmp(m,n,l) \label{eqn:Gas-BoundsCmp-TP} \\
  0 \leq f^{Cmp,H2}_{rp,k,m,n,l} \leq f^{Cmp,CH4}_{rp,k,m,n,l} \overline{B}^{H2} 
  \quad \forall rp,k,cmp(m,n,l) \label{eqn:H2-BoundsCmp-TP}
\end{align}
\end{subequations}

Finally, equations (\ref{eqn:Balances-INC}) represent the hydrogen (\ref{eqn:H2-Balance-INC}) and natural gas (\ref{eqn:CH4-Balance-INC}) balances.

\begin{subequations}
\label{eqn:Balances-INC}
\begin{align}
  \sum_{h2um(h2g,m)}    p^{H2}_{rp,k,h2g}
 +\sum_{h2um(h2p,m)}    p^{H2}_{rp,k,h2p}       
 +\sum_{h2um(h2s,m)}    p^{H2}_{rp,k,h2s}       \nonumber\\
 +\sum_{cl}             h2ns_{rp,k,m,cl}           
 +\sum_{mnl(n,m,l)}     f^{H2}_{rp,k,m,n,l} 
 -\sum_{mnl(m,n,l)}     f^{H2}_{rp,k,m,n,l}     \nonumber\\ 
 +\sum_{cmp(n,m,l)}     f^{Cmp,H2}_{rp,k,m,n,l} 
 -\sum_{cmp(m,n,l)}     f^{Cmp,H2}_{rp,k,m,n,l} 
 =
  \sum_{cl}             D^{H2}_{rp,k,m,cl}      \nonumber\\ 
 +\sum_{cl}             d^{Gas,H2}_{rp,k,m,cl}                           
 +\sum_{h2um(h2s,m)}    cs^{H2}_{rp,k,h2s}         
 +\sum_{gm(t=gas,m)}   (cs^{H2,E}_{rp,k,t}+cs^{H2,Aux}_{rp,k,t})      \nonumber \\
 +\sum_{cmp(m,n,l)}     CS^{Cmp,H2}_{m,n,l} f^{Cmp,H2}_{rp,k,m,n,l} 
 \quad \forall rp,k,m \label{eqn:H2-Balance-INC} \\
  \sum_{ch4um(ch4w,m)}  p^{CH4}_{rp,k,ch4w}
 +\sum_{ch4um(ch4s,m)}  p^{CH4}_{rp,k,ch4s} 
 +\sum_{mnl(n,m,l)}     f^{CH4}_{rp,k,m,n,l}        \nonumber\\
 -\sum_{mnl(m,n,l)}     f^{CH4}_{rp,k,m,n,l} 
 +\sum_{cmp(n,m,l)}     f^{Cmp,CH4}_{rp,k,m,n,l} 
 -\sum_{cmp(m,n,l)}     f^{Cmp,CH4}_{rp,k,m,n,l}    \nonumber\\
 +\sum_{cl}             ch4ns_{rp,k,m,cl}                              
 =
  \sum_{cl}             d^{CH4}_{rp,k,m,cl} 
 +\sum_{gm(t=gas,m)}   (cs^{CH4,E}_{rp,k,t}+cs^{CH4,Aux}_{rp,k,t}) \nonumber\\
 +\sum_{ch4um(ch4s,m)}  cs^{CH4}_{rp,k,ch4s}
 +\sum_{ch4um(h2p,m)}   cs^{CH4}_{rp,k,h2p}         \nonumber\\
 +\sum_{cmp(m,n,l)}     CS^{Cmp,CH4}_{m,n,l} f^{Cmp,CH4}_{rp,k,m,n,l}
 \quad \forall rp,k,m \label{eqn:CH4-Balance-INC}
\end{align}
\end{subequations}

\subsection{Power system policy}
\label{subsec:Policy}
Worldwide, power systems are transitioning towards high shares of renewables. Austria, for example, recently presented legislation \cite{EAG2021} setting the target of covering 100\% of national demand (on a net annual basis) from renewables by 2030. Motivated by this, we formulate a renewable power system constraint (\ref{eqn:RClean}) \cite{Wogrin2020} based on the minimum renewable generation rate $\kappa$. This constraint limits fossil-fired power generation and thus ensures that the specified renewable generation rate is met (or exceeded if it is optimal from a cost perspective). On the LHS, only the specific natural gas consumption associated to power generation of gas-fired units is considered. With this, supplying power demand from hydrogen-based thermal generation is permitted. We want to point out that (\ref{eqn:RClean}) is set up on an annual basis and therefore represents a net minimum. An exception to this is the 100\% case, where fossil-based generation is constrained to zero. 

In the framework, electricity demand from hydrogen units represents an additional (variable) demand on top of the baseline demand $D^E_{rp,k,i}$. The way constraint (\ref{eqn:RClean}) is designed implies that this variable demand can only be covered from variable renewable energy sources, battery energy storage systems (BESS), and hydrogen-based thermal generation.
However, it is most unlikely that hydrogen is produced based on electricity from the latter two technologies.
There are two reasons for this. First, BESSs entail losses, so it is more economical to produce hydrogen directly from renewable electricity. Second, it makes no sense to burn hydrogen for the sake of producing hydrogen -- again, due to losses.
\begin{align}
  \sum_{rp,k} W_{rp}^{RP} W_k^K 
                    \Big(   \sum_{t \neq gas} p^E_{rp,k,t}
                           +\sum_{t =    gas} cs^{CH4}_{rp,k,t} H^{CH4} / CS_{t}^{V}
                    \Big)
  \leq \nonumber\\ (1-\kappa) \sum_{rp,k,i} W_{rp}^{RP} W_k^K D^E_{rp,k,i} 
  \label{eqn:RClean}
\end{align}

\section{Case studies}
\label{sec:CaseStudies}
In this section, we apply the proposed ESM in two case studies and illustrate its potential for energy system planning and analyses.

First, we study the impact of the gas flow formulation (B-TP versus the novel B-PP) on generation expansion planning decisions in the power and hydrogen sectors, pipeline transmission expansion planning, and in an operational context. Our results indicate that the B-TP provides an acceptable approximation for generation expansion planning, while it lacks in terms of transmission expansion planning and cannot guarantee operational feasibility of the planned energy system.

In the second case study, we investigate the optimal ramp-up of the hydrogen sector on the path towards climate neutrality.
Especially during the early stages of this process, it is highly unclear how the hydrogen sector should optimally evolve, as neither a dedicated hydrogen demand nor production/storage infrastructure nor dedicated transmission infrastructure exists today.
Motivated by this, we utilize the flexible framework of the presented ESM to study the optimal deployment of hydrogen in the context of increasing levels of power and gas sector integration. This ranges from its deployment for the purpose of (long-term) energy storage in the power sector to a extensively sector-coupled energy system in which hydrogen can also be deployed as a substitute for natural gas in various economic sectors, thereby reducing CO\textsubscript{2} emissions.
Our results highlight the critical role of hydrogen transmission for ramping-up the hydrogen sector, showcase the impact of inter-sectoral effects in the power and gas sectors, and evaluate the effectiveness of CO\textsubscript{2} pricing on fostering the deployment of hydrogen in the gas sector.

\begin{figure}[!h]
\centerline{\includegraphics[scale=0.86]{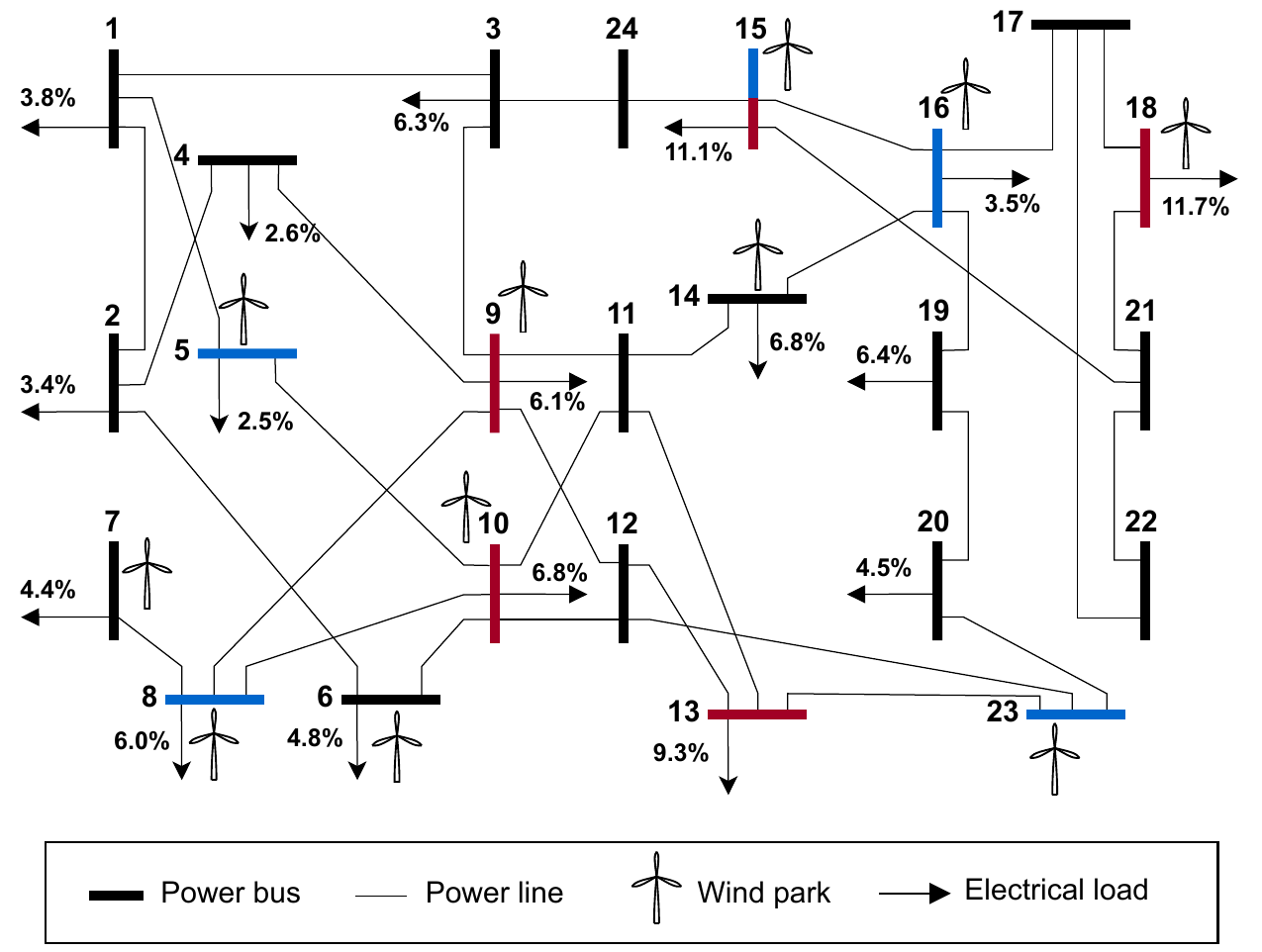}}
\caption{Graphical representation of the 23-bus power system.}
\label{fig:PowerSystem}
\end{figure}

All case studies are based on a modified version of an integrated 24-bus IEEE Reliability Test System and a 12-node gas system presented in Ordoudis et al. \cite{Ordoudis2019}. The energy system with a particular focus on the gas infrastructure is depicted in Fig.~\ref{fig:GasSystem}. Note that the model itself also represents the electric power system in detail, which is depicted in Fig.~\ref{fig:PowerSystem}, but since the original contributions of this paper lie within the formulation of the gas sector and its coupling with the power sector, we describe the gas sector in more detail. All case studies are solved on a notebook with a 2.80~GHz 11\textsuperscript{th} Generation Intel Core i7-1165G7 (4 cores) and 32 GB RAM using GAMS 37.1.0 and Gurobi 9.5.0. At the beginning of this section, we give the reader an overview of the most relevant model input data.

\subsection{Input data}
\label{subsec:Data}
This section provides an overview of the input data to the ESM at hand and the general setup for the case studies. Detailed input data is provided in an \href{https://github.com/tklatzer/ESM-data}{online appendix}.

The temporal framework for all case studies comprises seven representative days, which are determined by a k-medoids clustering procedure. The time series for power and natural gas demands are based on the Austrian system demands in 2020, scaled to the test system, and distributed to the buses and nodes respectively.
Power flows in the transmission system are represented by a DC-OPF approximation based on voltage angles.
\begin{figure}[!h]
\centerline{\includegraphics[scale=0.86]{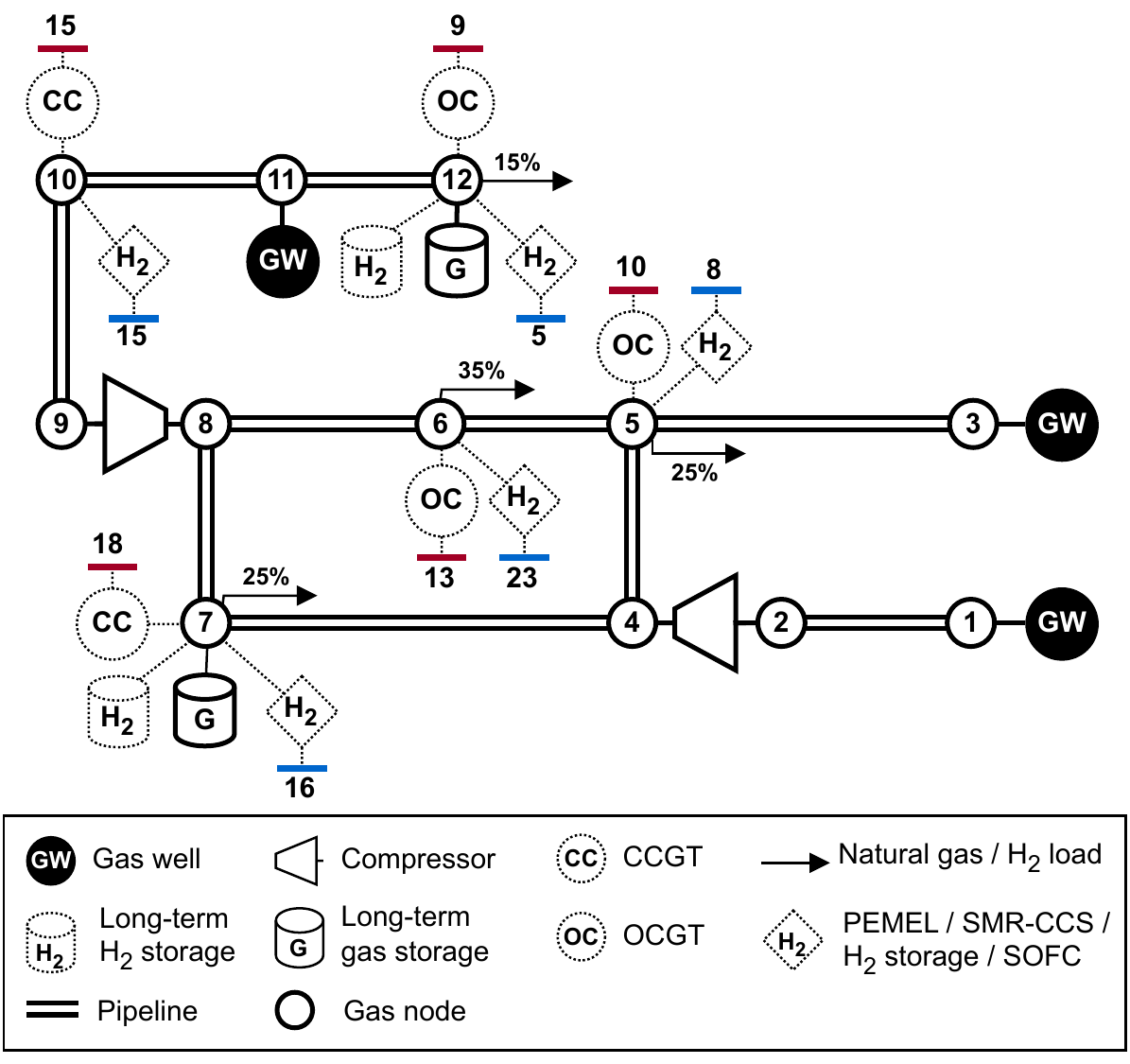}}
\caption{Graphical representation of the natural gas and hydrogen system.}
\label{fig:GasSystem}
\end{figure}
\newpage
The following sections present the considered ESM with the focus on the gas system.
We start from a brownfield setup\footnote{This is inspired by the European energy transition, which does not start from scratch either. However, the model also allows for a greenfield approach.} where three gas wells with maximum gas delivery rates of 0.425~MSm\textsuperscript{3}/h (nodes 1 and 3) and 0.283~MSm\textsuperscript{3}/h (node 11) supply the natural gas demand. In addition, two depleted natural gas fields (Table~\ref{tab:ParamStorage}) provide long-term natural gas storage.

The nodes of the gas system are connected by a high-pressure (43-68~bar) pipeline network. The network data for pipeline and compressor units is provided in Table~\ref{tab:ParamGasNet} and Table~\ref{tab:ParamComp} respectively.
In our previous work \cite{Klatzer2022}, we found that in existing expansion planning literature gas system parameters (and gas characteristics) tend to only be vaguely described, e.g. by Waymouth constants, or not at all.
Since the model (and the input data file) is available as open source, we try to provide as much information and flexibility to potential users (and ourselves) as possible, e.g. to implement real system data or to conduct sensitivity analyses.
Thus, the framework allows to specify all relevant pipeline and gas parameters for (steady-state) pipeline gas transmission in great detail.
This includes, inter alia, gas velocity, temperature, density, pipeline length, diameter, roughness, the Reynolds number (\ref{eqn:Reynolds}), the friction factor, etc. for each pipeline.
The maximum pipeline transmission capacity in Table~\ref{tab:ParamGasNet} is determined based on the steady-state gas flow equation (\ref{eqn:GasSteady}). The friction factor is typically described by the Colebrook-White equation, which is an implicit formulation and thus problematic for ESMs. However, the Chen equation (\ref{eqn:Chen}) gives a very good explicit approximation \cite{Domschke2017}.

The 12-node gas system at hand comprises two radial flow turbocompressor units with typical compression ratios $\Lambda^{sqr}$ of 1.20 and 1.30 respectively. The CUs consume approximately 0.2\% of the transported gas per 100~km to provide the required compression work.

\begin{subequations}
\begin{align}
Re = \frac{D v^M \rho^M}{\eta^M} \label{eqn:Reynolds}\\
\overline{F}^{Gas} = \sqrt{\bigg( \overline{P}_{m}^2 - \underline{P}_{n}^2 \bigg)\, \frac{1}{\lambda}\,\frac{D^5}{L}\,\frac{\pi^2}{16}\,\frac{T^N}{T^M}\,\frac{1}{p^N}\,\frac{1}{\rho^N}\,\frac{1}{K^M}} \label{eqn:GasSteady}\\
\frac{1}{\sqrt{\lambda}} = 
-2 \log_{10} \Bigg[\frac{\frac{\epsilon}{D}}{3.7065} - \frac{5.0425}{Re}
\log_{10} \bigg( \frac{\left(\frac{\epsilon}{D}\right)^{1.1098}}{2.8257} + \frac{5.8506}{Re^{0.8981}}\bigg) \Bigg]  \label{eqn:Chen}
\end{align}
\end{subequations}

\begin{table}[h]
\centering
\begin{tabular}{l|c|c|c|c|}
        & $L$     & $D$   & $R^{Gas}$                                                                                   & $\overline{F}^{Gas}$        \\
        & (km)    & (m)   & (10\textsuperscript{-5}((MSm\textsuperscript{3})\textsuperscript{2}/bar\textsuperscript{2}))& (MSm\textsuperscript{3}/h)  \\ \hline
 1--2   & 70      & 0.6   &  6.808                                                                                      & 0.435                       \\ \hline
 3--5   & 70      & 0.6   &  6.808                                                                                      & 0.435                       \\ \hline
 4--5   & 60      & 0.6   &  7.942                                                                                      & 0.469                       \\ \hline
 5--6   & 45      & 0.6   & 10.590                                                                                      & 0.542                       \\ \hline
 4--7   & 70      & 0.6   &  6.808                                                                                      & 0.435                       \\ \hline
 6--8   & 80      & 0.6   &  5.957                                                                                      & 0.407                       \\ \hline
 7--8   & 80      & 0.6   &  5.957                                                                                      & 0.407                       \\ \hline
 9--10  &125      & 0.6   &  3.812                                                                                      & 0.325                       \\ \hline
10--11  & 90      & 0.6   &  5.295                                                                                      & 0.383                       \\ \hline
11--12  & 85      & 0.6   &  5.606                                                                                      & 0.394                       \\ \hline
\end{tabular}
\caption{Parameters for pipelines.}
\label{tab:ParamGasNet}
\end{table}

\begin{table}[h]
\centering
\begin{tabular}{l|c|c|c|}
        & $\Lambda^{sqr}$ & $\overline{\Lambda}^{sqr}$  & $CS^{Cmp}$  \\
        & (p.u)           & (bar)                       & (p.u.)      \\ \hline
 2--4   & 1.20            & 30                          & 0.0015      \\ \hline 
 9--8   & 1.30            & 30                          & 0.0020      \\ \hline
\end{tabular}
\caption{Parameters for compressor units.}
\label{tab:ParamComp}
\end{table}

In the presented ESM, gas-fired thermal power plants interlink the gas and the power systems. Thermal candidate units (Table~\ref{tab:ParamThermal}) comprise highly efficient 400~MW CCGTs and flexible 200~MW OCGTs. We indicate the connection points of the thermal units (red) and hydrogen units (blue) to the power system by the numbers in Fig.~\ref{fig:GasSystem}.
We want to point out, that the investment options for CCGTs are limited to buses 15 and 18 as these are the buses with the highest power demand. Furthermore, the upper bound for the investment in thermal candidate units is one and the investment decision is binary. Finally, thermal units are subject to unit commitment (UC) and ramping constraints in all studied cases.

For the power sector, we consider 100~MW increments of solar and wind parks (Table~\ref{tab:ParamVRE}), and 50~MW increments of BESSs (Table~\ref{tab:ParamBESS}) as candidate units. The maximum hourly power generation potential of solar and wind units depends on the availability of solar and wind resources, which is express as capacity factors per bus and technology. The applied capacity factors reflect real time series for the Austrian system and are downloaded from \href{https://www.renewables.ninja/}{Renewables.ninja}.
For solar, we consider an expansion potential of 700~MW per bus. This is as solar resources are relatively evenly distributed and solar can be installed in various forms, e.g. on rooftops or as large-scale ground-mounted plants.
In contrast, wind parks typically have regional limitations, e.g. due to sufficient availability of wind, but also due to spatial and societal constrains. In order to reflect this, the expansion potential of wind is limited to busses 5-10, 14-16, 18, and 23 and 600~MW per bus.
With the transition of power systems to a high share of renewables, storage technologies are becoming increasingly important for the operation of the system. To this end, we assume an extensive expansion potential of 750~MW of BESSs per bus.

\begin{landscape}
\begin{table}[h]
\centering
\begin{tabular}{l|c|c|c|c|c|c|c|c|c|}
           & $\underline{P}^E$  & $\overline{P}^E$  & $RU,RD$       & $E$                        & $CS^{V}$  & $CS^{SU}$& $CS^{UP}$ & $C^{INV}$         \\
           & (MW)               & (MW)              & (MW)          & (tCO\textsubscript{2}/MWh) & (p.u.)    & (GWh)    & (GWh/h)   & (\euro/MW/y)      \\ \hline
CCGT       &  80                &  400              &  $\pm$ 160    &  0.181                     & 2.092     & 1.162    & 0.349     & 41,819            \\ \hline
OCGT       &  20                &  200              &  $\pm$ 180    &  0.181                     & 2.324     & --       & 0.166     & 24,781            \\ \hline
\end{tabular}
\caption{Parameters for gas-fired generation units.}
\label{tab:ParamThermal}
\end{table}

\begin{table}[!h]
\centering
\begin{tabular}{l|c|c|c|c|}
           & $\overline{P}^E$   & $C^{OM}$      & $C^{INV}$         \\
           & (MW)               & (\euro/MWh)   & (\euro/MW/y)      \\ \hline
Wind       & 100                & 2.0           & 72,642            \\ \hline
Solar      & 100                & --            & 84,467            \\ \hline
\end{tabular}\caption{Parameters for wind and solar units.}
\label{tab:ParamVRE}
\end{table}

\begin{table}[!h]
\centering
\begin{tabular}{l|c|c|c|c|c|c|c|c|c|c|}
        & $\overline{P}^{E}$   & $\eta^{CH/DIS}$   & $C^{OM}$       & $C^{INV,Pow}$ & $C^{INV,En}$  & $ETP$  \\
        & (MW)                 & (p.u.)            & (\euro/MWh)    & (\euro/MW/y)  & (\euro/MWh/y) & (h)    \\ \hline
BESS    & 50                   & 0.922             & 4.0            & 56,667        & 13,333       & 4      \\ \hline
\end{tabular}
\caption{Parameters for battery energy storage systems.}
\label{tab:ParamBESS}
\end{table}
\end{landscape}

For hydrogen production, we consider two fundamentally different production principles: water electrolysis and natural gas reforming. 
As of today, alkaline electrolysis (AEL) is the most common electrolysis-based hydrogen production technology. However, proton exchange membrane electrolyzer (PEMEL) units (Table~\ref{tab:ParamEL}) are expected to be the predominant electrolysis-based technology in the near future \cite{Bohm2020}. Although CAPEX for PEMEL is still higher than for AEL, a decisive technological advantage of PEMEL is its startup and ramping characteristics, which enables fast load changes, e.g. to adapt to the generation pattern of renewables.
As pointed out in section \ref{subsec:H2-sector}, 75\% of today's hydrogen demand is supplied from steam-methane reforming, which is a carbon-intensive process. To mitigate the bulk of carbon emissions, we consider steam-methane reforming units with carbon capture and storage (SMR-CCS) (Table~\ref{tab:ParamSMR}), which is already at technology readiness level (TRL) nine, according to Fan et al. \cite{Fan2021}.
Data for the 20~MW PEMEL and the 50,000 Sm\textsuperscript{3}/h SMR-CCS units\footnote{The capacity of SMR-CCS is scaled down to the test system, assuming a linear relationship of capacity to investment costs.} is based on \cite{IRENA2018} and \cite{IEAGHG2017} respectively.

A key characteristic of hydrogen is its storage potential. In the present ESM, we consider high-pressure steel tanks for short-term and salt caverns for long-term, e.g. inter-seasonal, hydrogen storage (Table~\ref{tab:ParamStorage}). The investment costs per installed power and energy capacity of the two technologies under consideration are determined with respect to the maximum consumption rate.
Data for steel tanks and salt caverns is based on \cite{Gorre2018} and \cite{Gorre2018,Kruck2013} respectively.

Finally, solid oxide fuel cells (SOFCs), currently at TRL 6-7 \cite{Schmidt2020}, complete the range of considered hydrogen technologies. Their high operating temperature (700-1000°C) offers additional potential for sector coupling, e.g. for combined heat and power. Data is provided in Table~\ref{tab:ParamFC} and based on \cite{Schmidt2020}.

For the purpose of a comprehensive hydrogen investment portfolio, we consider 400~MW of PEMEL, 500,000~Sm\textsuperscript{3}/h of SMR-CCS, 175,000~Sm\textsuperscript{3}/h of steel tanks (based on the maximum consumption), and 66,000~Sm\textsuperscript{3}/h of SOFCs at each node of the gas network. Salt caverns, on the other hand, are subject to geological requirements and thus, candidate units are limited to nodes 7 and 12. Just as with gas-fired power plants, the upper limit for investments in salt caverns is one and the investment decision is binary.
In contrast to the binary investment decisions for thermals and salt caverns, investments in wind, solar, BESSs, PEMEL, SMR-CCS, hydrogen steel tanks, and SOFCs are continuous, which is a good approximation since we are planning a GW-scale energy system. 

\begin{table}[!h]
\centering
\begin{tabular}{l|c|c|c|c|}
           & $\overline{P}^E$   & $HPE$                         & $C^{OM}$                          & $C^{INV}$         \\
           & (MW)               & (Sm\textsuperscript{3}/MWh)   & (\% of C\textsuperscript{INV})    & (\euro/MW/y)      \\ \hline
PEMEL      &  20                & 213.91                        & 2.0                               & 35,000            \\ \hline
\end{tabular}
\caption{Parameters for electrolyzer units.}
\label{tab:ParamEL}
\end{table}

\begin{table}[h]
\centering
\begin{tabular}{l|c|c|c|c|c|}
           & $\overline{P}^{H2}$        & $HPC$    & $E$                                                                & $C^{OM}$                       & $C^{INV}$                                \\
           & (Sm\textsuperscript{3}/h)  & (p.u.)   & (kgCO\textsubscript{2}/Sm\textsuperscript{3}H\textsubscript{2})    & (\% of C\textsuperscript{INV}) & (\euro/(Sm\textsuperscript{3}/h)/y)      \\ \hline
SMR-CCS    & 50,000                     & 0.69     & 0.09                                                               & 2.9                            & 159.39                                   \\ \hline
\end{tabular}
\caption{Parameters for steam-methane reforming units with carbon capture and storage.}
\label{tab:ParamSMR}
\end{table}

\begin{table}[!h]
\centering
\begin{tabular}{l|c|c|c|c|}
           & $\overline{CS}^{H2}$       & $EPH$                         & $C^{OM}$                          & $C^{INV}$                             \\
           & (Sm\textsuperscript{3}/h)  & (kWh/Sm\textsuperscript{3})   & (\% of C\textsuperscript{INV})    & (\euro/(Sm\textsuperscript{3}/h)/y)   \\ \hline
SOFC       &  3,300                     & 1.797                         & 2.0                               & 699                                   \\ \hline
\end{tabular}
\caption{Parameters for fuel cell units.}
\label{tab:ParamFC}
\end{table}

\begin{landscape}
\begin{table}[h]
\centering
\begin{tabular}{l|c|c|c|c|c|c|c|c|c|}
                     & $\overline{P}^{CH4,H2}$    & $\overline{CS}^{CH4,H2}$   & $\eta^{CH/DIS}$   & $InRes^{H2}$ & $\underline{R}^{H2}$  & $C^{OM}$                       & $C^{INV,Pow}$                          & $C^{INV,En}$                       & $ETP$  \\
                     & (Sm\textsuperscript{3}/h)  & (Sm\textsuperscript{3}/h)  & (p.u.)            & (p.u.)       & (p.u.)                & (\% of C\textsuperscript{INV}) & (M\euro/(MSm\textsuperscript{3}/h)/y)  & (M\euro/MSm\textsuperscript{3}/y)  & (h)    \\ \hline
Depleted NG field    & 250,000                    & 180,000                    & 0.995             & 0.80         & 0.60                  & --                             & --                                     & --                                 & 500    \\ \hline
Hydr. salt cavern    & 130,000                    & 130,000                    & 0.995             & 0.78         & 0.55                  & 2.0                            & 2.50                                   & 1.87                               & 362    \\ \hline
Hydr. steel tank     &   5,000                    &   3,500                    & 0.995             & --           & --                    & 1.5                            & 3.75                                   & 1.25                               & 12     \\ \hline
\end{tabular}
\caption{Parameters for natural gas and hydrogen storage units.}
\label{tab:ParamStorage}
\end{table}
\end{landscape}

\clearpage
\subsection{Impact of the gas flow formulation on planning results}
\label{subsec:Quality}
As mentioned in the introduction, the EU has ambitious plans for the expansion of the hydrogen sector \cite{RePowerEU2022,RePowerEUPlan2022}. During the early stages of this process, utilizing the gas system to transport hydrogen via blending will be key. Against this background, this case study is focused on assessing the impact of the gas flow formulation on planning results in the integrated sector-coupled energy system. In particular, we utilize the flexible structure of the model and compare the B-TP, which omits the physical relation between gas flow and pressure, versus the novel B-PP formulation and their implications on generation expansion planning (GEP) decisions in the hydrogen and power sectors, pipeline transmission expansion planning (TEP), as well as system operation. The quality of the planning results is evaluated in terms of the regret (measured in the form of non-supplied hydrogen) that results from fixing all investments from the B-TP framework and re-running the model with the gas flows governed by the more realistic B-PP. In the following, we describe the specific assumptions made for this case study.

As a basis for modeling the hydrogen sector, we consider the future hydrogen demand, e.g. of the iron \& steel and the chemical industries, as an exogenous parameter. Since no time series for large-scale hydrogen demand are publicly available to date, we assume that the hydrogen demand follows the same time series and local distribution as the natural gas demand scaled down to the test system. 

To reflect the decarbonization of the power sector, we set the policy requirement that at least 95\% of the total generated electricity must originate from renewable sources or, in other words, thermal generation is limited to at most 5\% of the total generation. That way, we can ensure that the produced hydrogen via electrolysis qualifies as renewable, at least under the currently available EU draft delegated regulation \cite{DefRenewableH2Draft2022}. At the same time, this still permits the operation of natural gas-fired units, which ensures that their operational characteristics, e.g. startup, ramping etc., and the resulting consumption from the gas system are captured and accounted for in this case study.

The GEP candidate units in the power and hydrogen sectors coincide with the portfolio described in section \ref{subsec:Data}. The pipeline connecting nodes 5 and 6 of the gas transmission system (see Table~\ref{tab:ParamGasNet}) represents a candidate pipeline for TEP capable of natural gas and hydrogen blending. The associated investment cost is 27~M\euro\footnote{This is in line with the benchmark costs for pan-European natural gas transportation presented in \cite{CEER2019b}.} (binary investment decision), which is annualized based on an annuity factor of 5\%. Finally, the MILP gap is set to 1\%.

We start this case study with an expansion planning problem (GEP and TEP) with the gas flows governed by the B-TP versus the B-PP formulation. 
As pointed out in section~\ref{subsec:GasNetwork}, gas flows under the B-PP framework are non-linear and non-convex and therefore have to be linearized. For the piecewise linearization (see Fig.~\ref{fig:B-PP}) we consider a total of 6 increments.
Finally, the maximum permitted hydrogen blending rate is 10\% of the actual natural gas flow. 

\begin{figure}[t]
\centerline{\includegraphics[scale=0.45]{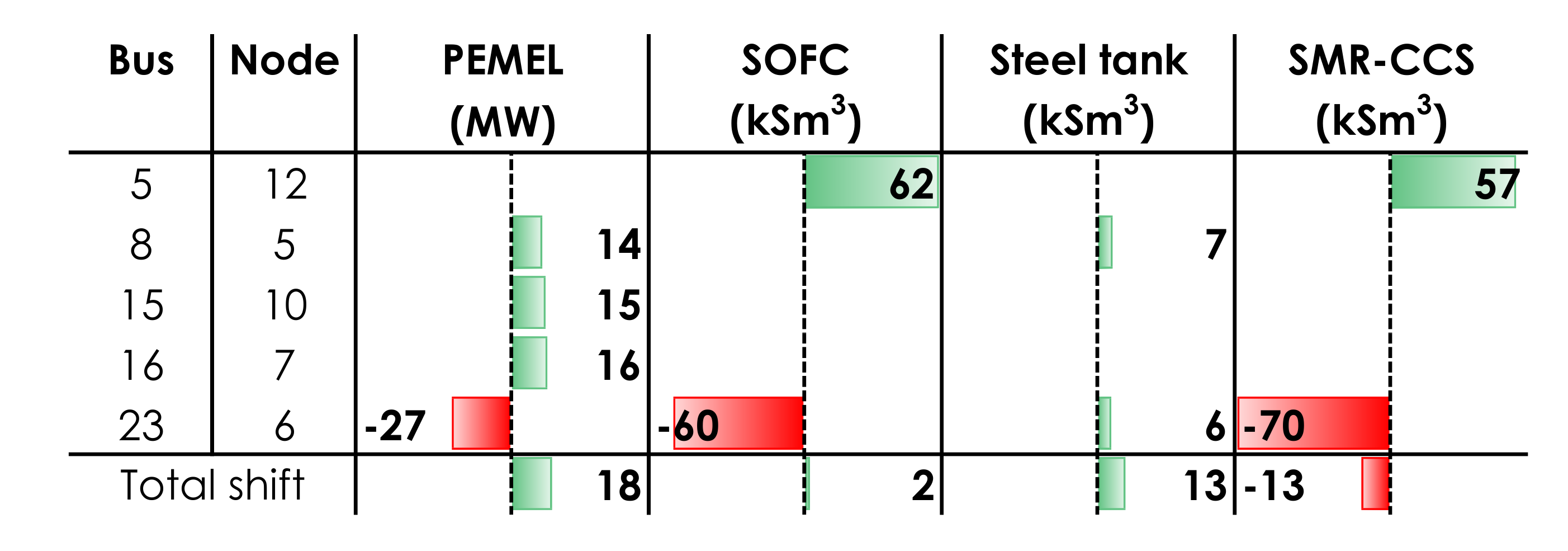}}
\caption{Difference in installed hydrogen capacity per bus and node under the B-PP compared to the B-TP.}
\label{fig:SitingH2}
\end{figure}

Expansion planning under the two frameworks results in total system costs of 1,094~M\euro\hspace{1pt} for the B-TP and 1,107~M\euro\hspace{1pt} for the B-PP.
Given the relatively small difference, it appears that the B-TP performs quite well within a GEP context.
However, taking a closer look at the operational results reveals significant changes when planning under the more realistic B-PP compared to the B-TP framework, e.g. the reversal of gas flows in pipelines, which is a direct consequence of linking gas flows with gas pressure. For the energy system at hand, this results in a shift of total natural gas production from gas well 3 (2010$\rightarrow$1469~MSm\textsuperscript{3}) to gas wells 11 (437$\rightarrow$919~MSm\textsuperscript{3}) and 1 (180$\rightarrow$214~MSm\textsuperscript{3}) under the B-PP.
This in turn affects siting decisions of hydrogen infrastructure as the capacity to transport hydrogen in a pipeline via blending is a function of the actual natural gas flow. Fig.~\ref{fig:SitingH2} depicts the difference in installed hydrogen capacity under the B-PP compared to the B-TP framework.
In particular, SMR-CCS (and SOFC) capacities are quasi completely relocated from node~6 to node~12 (this is in line with the shift of gas production), while PEMEL capacities are distributed across multiple nodes.
Ultimately, it is the relocation of PEMEL capacities that is the main reason behind the increased total system cost when planning under the B-PP. The costs, however, are not borne in the hydrogen but in the power sector.
This is because powering PEMELs (and the power sector) requires the expansion of renewables at sites with less favourable capacity factors compared to the B-TP, at least in the studied system.
As a general takeaway, planning under the B-PP can lead to a shift from the expansion of centralized (SMR-CCS) towards decentralized hydrogen production units (PEMEL). 

\begin{figure}[!b]
\centerline{\includegraphics[scale=0.40]{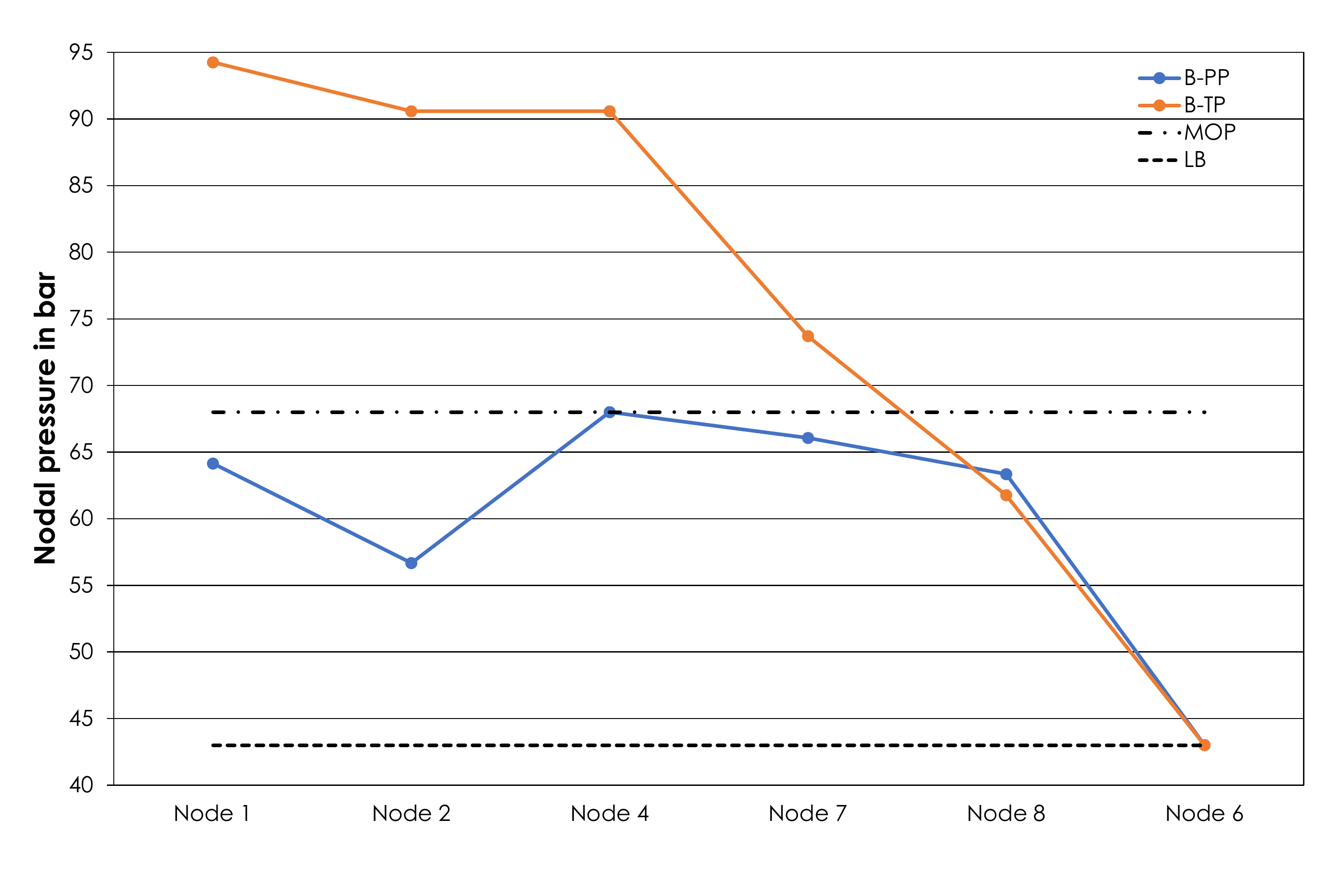}}
\caption{Nodal pressures on the pipeline stretch from node~1 to node~6 determined under the B-PP versus the B-TP framework. LB corresponds to the lower pressure bound and MOP to the maximum operating pressure.}
\label{fig:Pressure}
\end{figure}

When it comes to TEP, the investment decision for the specified candidate pipeline changes under the two frameworks (0 for the B-TP; 1 for the B-PP)\footnote{For a system-wide blending rate of zero there is no TEP under both frameworks. This indicates that hydrogen blending per se can be sufficient to trigger investments in pipeline infrastructure, at least for the system at hand.}, which is significant.
Against this background, we assess the quality of the planning decisions made on the basis of regret, which we quantify in terms of non-supplied hydrogen. To this end, we fix the optimal investments (GEP and TEP) determined under the B-TP framework and run an operational problem with the gas flows governed by the more realistic B-PP, which results in a total of 23~MSm\textsuperscript{3} of non-supplied hydrogen\footnote{For the cost of non-supplied hydrogen, we assume 3 \euro/Sm\textsuperscript{3}.} (4\% of total hydrogen deployment) and total system costs of 1,162~M\euro\hspace{1pt}.
The reason for this is that the gas flows determined under the the B-TP framework violate the maximum operating pressure (MOP) of the gas transmission system to a large extent.
Fig.~\ref{fig:Pressure} depicts the nodal pressures which would result from the gas flows determined under the B-TP framework along the pipeline stretch connecting node~1 (gas well~1) and node~6 (highest natural gas and hydrogen demand).
From these results it is clear that the gas flows derived under the B-TP are not suitable for operational planning as fixing them would result in an infeasible system.

Nevertheless, a strong argument in favour of using the B-TP over the B-PP is its much lower computational burden as reflected in the number of variables (75,465 (2,598 discrete) for the B-TP and 101,001 (12,678 discrete) for the B-PP). The resulting computational time is highly case-dependent, but in general for the energy system at hand and a MILP gap of 1\% it is under 300 seconds for the B-TP and about ten hours for the B-PP.

Considering the above, we conclude that the B-TP is an acceptable approximation for the purpose of generation expansion planning and small hydrogen blending rates. However, when it comes to pipeline expansion planning and operational feasibility, model results can vary significantly depending on the case, and the B-PP can be a valuable tool for assessing the differences.

\subsection{Ramping up the hydrogen sector}
\label{subsec:RampingH2}
The motivation for this case study is to investigate the optimal ramp-up of the hydrogen sector in a sector-coupled power and gas system on the path towards climate neutrality.
The underlying idea is to leverage existing infrastructure, such as transmission pipelines for blending, and to deploy hydrogen in the power sector or gas sector or both, depending on where it has the greatest impact (in terms of total system cost). 
To this end, we utilize the flexible formulation of the proposed ESM to allow for the optimal substitution of hydrogen for natural gas based on maximum substitution rates per gas demand sector, rather than specifying a dedicated hydrogen demand for which no time series are available as of today.
Furthermore, we limit hydrogen production to PEMEL (and omit SMR-CCS) and we impose the green power system constraint (\ref{eqn:RClean}) with $\kappa=100\%$. This way we ensure that both the power sector and thus the produced hydrogen are totally renewable\footnote{The legal framework defining when hydrogen qualifies as renewable is still missing. To date, only a draft delegated regulation exists \cite{DefRenewableH2Draft2022}.}. 
The cost of natural gas is assumed as 0.097~\euro/Sm\textsuperscript{3}\footnote{This is in the range of TTF spot and year-ahead prices for natural gas in 2020 \cite{GasCoalCost2022}, which we assume for modeling the energy system in 2030.}.
Finally, for this case study, the MILP gap is set to 0.1\%. Besides that, the setup for this case study coincides with the general setup described in section~\ref{subsec:Data}.

This the case study is divided into three incremental sections. To give the reader a better understanding, the sections are structured in a way that the scope of hydrogen, and thus the integration of the energy system, increases. 
In section~\ref{sec:H2PowerSec}, we examine the deployment of hydrogen only in the power sector when hydrogen is limited to power-to-power technology and highlight the crucial role of hydrogen transmission.
In section~\ref{sec:H2PowerGas}, we study the ramp-up of the hydrogen sector in the context of a sector-coupled power and gas system, where we identify weak inter-sectoral effects from the gas sector towards the power sector, while these effects are more pronounced vice versa.
In the final section of this case study (\ref{sec:H2andCO2Emissions}), we apply the model to examine the impact of CO\textsubscript{2} pricing policies on the deployment of hydrogen in the gas sector, and find that the spatio-sectoral distribution of gas demand of the underlying energy system crucially affects its effectiveness.

\subsubsection{Hydrogen in the power sector}
\label{sec:H2PowerSec}
\begin{figure}[!h]
\centerline{\includegraphics[scale=0.40]{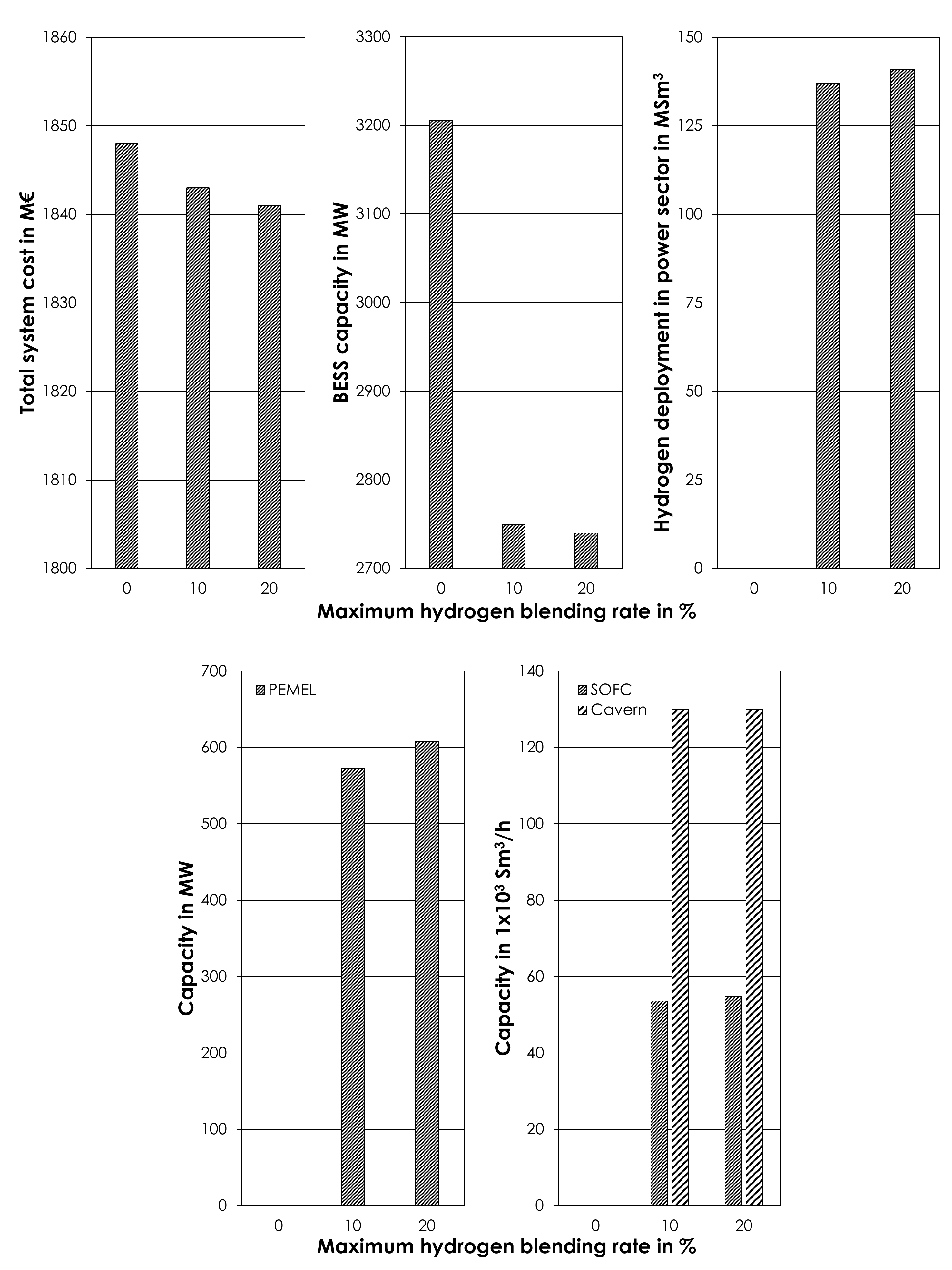}}
\caption{Ramp-up of the hydrogen sector for different maximum blending rates.}
\label{fig:H2PowerSector}
\end{figure}

Besides the massive expansion of renewable energies, storage technologies are an important cornerstone for the decarbonization of the power sector. Hydrogen represents a promising storage technology, as it can be produced from (renewable) electricity via electrolysis, transported, stored for days, months, or even seasons, and used to generate electricity on demand. Hence, in this section we explore the impact of deploying hydrogen in the power sector and study the impact of hydrogen transmission via blending on the ramp-up of the hydrogen sector. 
The results are summarized in Fig.~\ref{fig:H2PowerSector} and include, inter alia, the total system cost, the installed BESS capacity, the optimal hydrogen deployment in the power sector, etc. for different maximum blending rates (x-axis). 

We start with a case where the production, storage, and consumption of hydrogen are local.
Total system costs for this case comprise 1,848~M\euro. Since the power system is totally renewable ($\kappa\!=\!100\%$), curtailment of renewables increases to 101\% of total power demand. However, there are no investments in hydrogen infrastructure at all. This is since utilizing hydrogen as power-to-power technology requires investing in hydrogen technology chains consisting of EL, FC, and hydrogen storage for temporal shift of hydrogen, e.g. within a day or up to seasons. Ultimately, the hydrogen technology chain is competing against other storage technologies, e.g. BESS, which are preferred by the model (the more effective technology) based on the techno-economic assumptions for 4-hour lithium-ion BESS in 2030 in this case study \cite{Cole2021}. This is reflected by the investment in BESS, which is 3,206~MW (or 16.6\% of total installed capacity in the power system).
However, a potential advantage of hydrogen over BESS is its capacity to be transported.
In order to consider the spatial dimension, we enable hydrogen and natural gas blending for transmission in existing pipeline infrastructure, which is, for the sake of simplicity, governed by the B-TP formulation. We find that a blending rate of 10\%, which is currently permitted, e.g. in Austria~\cite{GasConnect2021}, is sufficient to trigger investments in hydrogen infrastructure.
The decisive factor for this is that blending facilitates the investment in the salt cavern at node 7 for long-term hydrogen storage.
As a result, total system costs decrease by 5~M\euro\ (-0.3\%), with BESS investments dropping by 456~MW (-14.2\%). Instead, the model invests in 573~MW PEMELs that produce 137~MSm\textsuperscript{3} hydrogen, of which 38~MSm\textsuperscript{3} (27.7\%) is transported via pipelines. Ultimately, hydrogen is utilized to generate 225~GWh of electricity, representing 19.8\% of total electricity generation from storage technologies.

For a blending rate of 20\% (currently under discussion in the UK \cite{UKH2Strategy2021}), total system costs decrease by another 2~M\euro\ (mainly to a shift of the installed capacity of PEMELs across nodes) and the amount of transported hydrogen increases to 116~MSm\textsuperscript{3} (82.3\% of the total hydrogen produced). 
However, our results show that the effectiveness of blending (in terms of total system cost) decreases significantly above 20\% and the deployment of hydrogen stagnates at 141~MSm\textsuperscript{3}, at least for the system at hand.
Moreover, it is noteworthy that the installed PEMELs only achieve 800--1,650 full-load hours (FLHs), averaging at 1,150 hours.

Finally, the attainable reduction of system costs through the deployment of hydrogen is not sufficient to trigger investments in dedicated hydrogen pipelines\footnote{This is for a case where the blending rate in the existing pipelines is zero, but the model has the option to invest in a dedicated hydrogen network running in parallel to the existing network (binary investment decision per pipeline, pipe diameter 0.35~m). The investment costs correspond to benchmark costs for pan-European natural gas transmission presented in \cite{CEER2019b} and are annualized on the basis of an expected lifetime of 60 years.}.

From the results above, we conclude that hydrogen transmission via blending can act as a crucial lever to initiate its deployment in the power sector, especially during the early stages of ramping up the hydrogen sector. However, as hydrogen competes against other storage technologies, its large-scale deployment solely as a power-to-power technology does not appear to be economic, as indicated by the limited number of FLHs.

\subsubsection{Hydrogen in the sector-coupled power and gas system}
\label{sec:H2PowerGas}
In addition to its application as a storage technology in the power sector, hydrogen can also be deployed to decarbonize the gas sector. Therefore, in the second part of this case study, we investigate the ramp-up of the hydrogen sector in a sector-coupled power and gas system and study inter-sectoral effects. 

To this end, and in addition to blending, we activate the model option to substitute a share of the natural gas demand with hydrogen, assuming that the resulting blend is combusted. Thus, the energy content of the blend has to be equivalent to the lower heating value of natural gas (see~(\ref{eqn:Gas-DemBlend})).
The flexible framework of the model allows the specification of different maximum volumetric hydrogen substitution rates for various sector classes. This is relevant since, e.g., household appliances are likely to be more limited in terms of hydrogen substitution than, e.g., processes in the industry sector, where higher substitution rates could be achieved more rapidly.\footnote{For the purposes of this case study, we segment the natural gas demand and assign it to the iron \& steel and chemical industry sectors (maximum substitution rate of 20\%), and others (10\%). Furthermore, the maximum blending rate for pipeline transmission is set to 20\% for this case study.}

First, we assume that the deployment of hydrogen is limited to the gas sector (by omitting FCs). Here we find that the underlying cost of natural gas (0.097~\euro/Sm\textsuperscript{3}) is on the verge of triggering investments in hydrogen infrastructure. A twenty-five percent increase of the natural gas cost results in the deployment of 110~MSm\textsuperscript{3}, and a fifty percent increase results in 140~MSm\textsuperscript{3} of hydrogen.
In the gas sector, hydrogen can be deployed continuously as a substitute for natural gas, which is reflected by the high number of FLHs of PEMELs ($\geq$~7,600 on average). 
Therefore, hydrogen storage plays only a subordinate role and thus the model only invests in short-term storage via steel tanks but no long-term storage.
This contrasts with the power-to-power case (see section \ref{sec:H2PowerSec}), where the model exclusively invests in a salt cavern for long-term hydrogen storage.

In the following, we consider the most complex case in which hydrogen can be deployed in the power and gas sectors (by including FCs). In this case, natural gas costs of 0.097~\euro/Sm\textsuperscript{3} are sufficient to trigger the deployment of hydrogen in the gas sector (154~MSm\textsuperscript{3}). However, interestingly, it is the power sector that drives the expansion of hydrogen infrastructure. This is evident from the investment decisions and their distribution, which is highly similar to the power-to-power case in section~\ref{sec:H2PowerSec} and the number of FLHs of PEMELs (2,200 hours on average).
In the sector-coupled case, increasing the cost of natural gas only results in a modest additional expansion of PEMEL capacity. However, the baseline deployment of hydrogen in the power sector remains unaffected by this and is continuously at 127~MSm\textsuperscript{3}. Thus, we conclude that the inter-sectoral effect of the gas sector on the power sector is weak with respect to hydrogen.
Vice versa, the inter-sectoral effect is more pronounced. For example, if the lifetime of lithium-ion BESS increases from 15 to 20 years for the same investment cost, not only does deploying hydrogen in the power sector become negligible ($\leq$ 1~MSm\textsuperscript{3}), but also its deployment in the gas sector decreases significantly (-37.0\%).

In our view, it is very likely that private companies will be at the forefront of the initial investments in hydrogen production infrastructure, as they can also establish demand in parallel. Moreover, as we have shown in our case study, the inter-sectoral effect from the gas towards the power sector is weak. Given the above, it is very likely that initially, the topology of hydrogen infrastructure will evolve similar to the gas sector-only case, which leads to the highest number of FLHs of PEMELs (see above).
In order to stimulate the topology of the future hydrogen system to evolve in the sense of cost-optimal energy system planning, it is likely that appropriate steering measures will be needed. As we have shown, such a system topology can ultimately foster the holistic ramp-up of the hydrogen sector and its deployment in both the gas and the power sectors.

\subsubsection{Impact of CO\textsubscript{2} pricing}
\label{sec:H2andCO2Emissions}
In the final part of this case study we analyze the impact of CO\textsubscript{2} pricing on the deployment of hydrogen in the gas sector.
As indicated above, the cost of natural gas can provide an incentive for the ramp-up of the hydrogen sector. However, it is the result of a global market and therefore difficult to estimate, regulate, and control (at least without policy intervention). Compared to that, CO\textsubscript{2} pricing represents a lever that can be applied in a more controlled, targeted, and predictive way, which is important for planning certainty.

Following this idea, we activate the model option to consider the cost for CO\textsubscript{2} emissions, which introduces an additional incentive to substitute natural gas with hydrogen.
The underlying emission reduction potential is based on average CO\textsubscript{2} emissions resulting from the combustion of natural gas in Austria \cite{ATNationalInventoryReport2022}.
During phase 4 (2021-2030) \cite{EUETSPhase42021} of the EU emission trading system, some industries (e.g. the iron \& steel industry and the chemical industry) will still receive 100\% of their determined emissions as free allowances, as they are considered industries at risk of carbon leakage \cite{SectorsAtRiskOfCarbonLeakage2019}. Thus, we exclude these industries, which yields average CO\textsubscript{2} emissions of 1.96~ktCO\textsubscript{2}/MSm\textsuperscript{3}CH\textsubscript{4} for the combustion of natural gas. However, due to the system topology, this introduces a dependency on the spatio-sectoral distribution of gas demand which can affect the ramp-up of the hydrogen sector.

\begin{figure}[h]
\centerline{\includegraphics[scale=0.40]{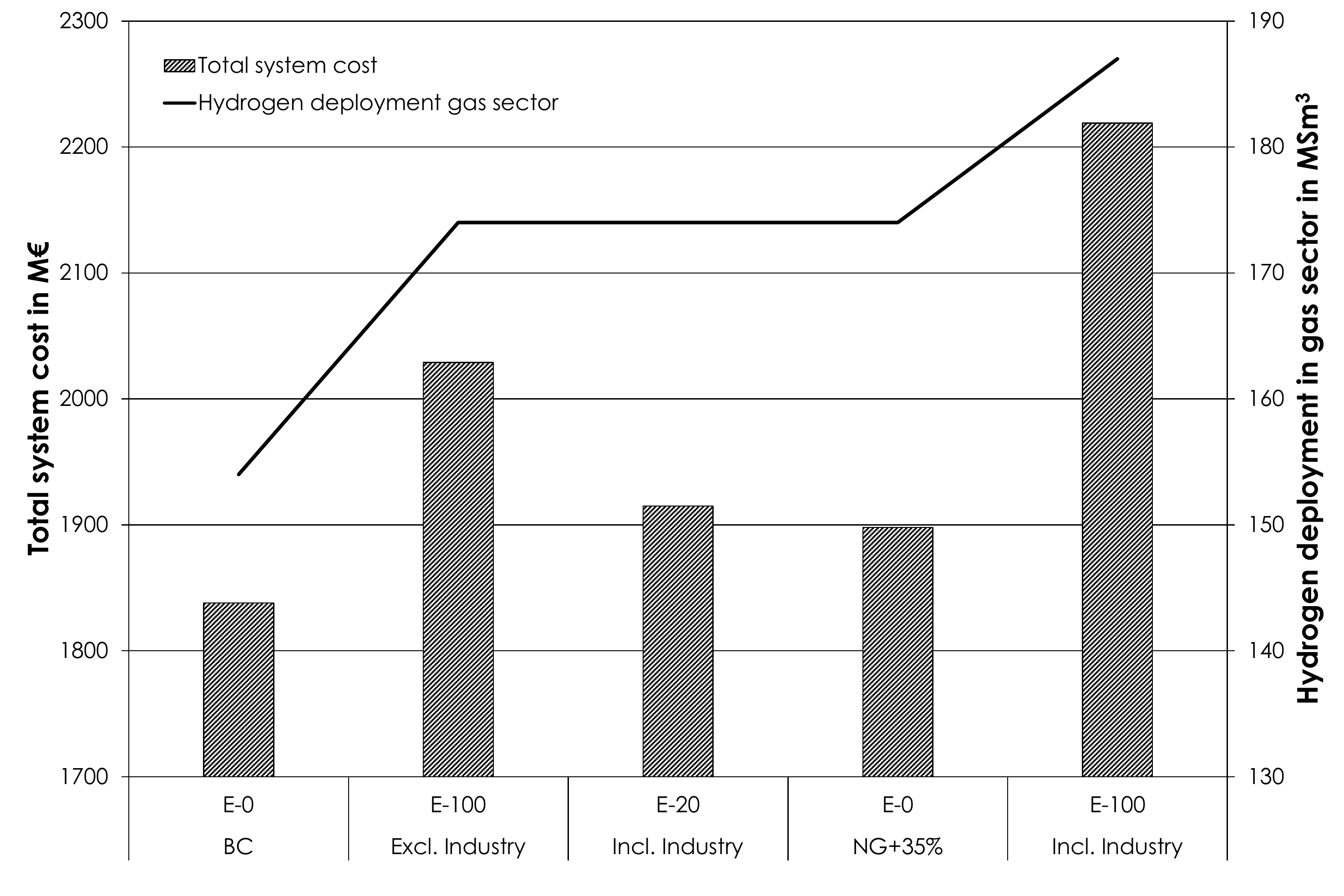}}
\caption{Deployment of hydrogen in the gas sector based on the cost of CO\textsubscript{2} allowances (E), the affected sectors, and the cost of natural gas (NG).}
\label{fig:RampupCO2}
\end{figure}

In the following, we quantify the impact of this spatial effect for the energy system at hand and analyze the effectiveness of CO\textsubscript{2} pricing to foster the deployment of hydrogen in the gas sector.
The results are depicted in Fig.~\ref{fig:RampupCO2}, which illustrates the deployment of hydrogen in the gas sector and the according total system cost as a function of the cost of CO\textsubscript{2} allowances (E).
Starting from an energy system where 154~MSm\textsuperscript{3} of hydrogen are deployed in the gas sector (base case (BC), E-0), we assume that the policy maker strives to stimulate the additional deployment of 20~MSm\textsuperscript{3} of hydrogen through CO\textsubscript{2} pricing. With the iron \& steel and chemical industries excluded from CO\textsubscript{2} pricing (Excl. Industry), achieving this policy goal requires a cost of 100~\euro/tCO\textsubscript{2} (E-100), whereas a cost of 20~\euro/tCO\textsubscript{2} is sufficient if both industries are included\footnote{For this case study, we assume average CO\textsubscript{2} emissions of 2.29~ktCO\textsubscript{2}/MSm\textsuperscript{3}CH\textsubscript{4} for the iron \& steel industry and 2.17~ktCO\textsubscript{2}/MSm\textsuperscript{3}CH\textsubscript{4} for the chemical industry based on \cite{ATNationalInventoryReport2022}.} (Incl. Industry, E-20). In terms of total system costs, this results in a difference of 114~M\euro, which reflects the costs of the spatio-sectoral distribution of gas demand inherent to the energy system at hand.
In contrast, the policy target could also be achieved by a system-wide increase of the cost of natural gas\footnote{The base cost of natural gas is assumed as 0.097~\euro/Sm\textsuperscript{3}.} by 35\% (NG+35\%), e.g., via taxation. In this case, the total system cost decreases by 17~M\euro, since this mechanism is independent of the spatio-sectoral distribution of gas demand of the underlying energy system.

Finally, and for the sake of completeness, if the policy maker establishes CO\textsubscript{2} costs of 100~\euro/tCO\textsubscript{2} for the entire gas sector (Incl. Industry, E-100), this would stimulate the deployment of 33~MSm\textsuperscript{3} of hydrogen.

From the above, we conclude that the spatio-sectoral distribution of gas demand of an energy system can have an significant impact on the effectiveness of CO\textsubscript{2} pricing policy measures.

\section{Conclusions}
\label{sec:Conclusions}
In this paper, we have presented a novel gas flow formulation for modeling detailed natural gas and hydrogen blending for pipeline gas transmission. The proposed blending pressure problem is based on a piecewise linearization of the non-linear and non-convex steady-state gas flow equation, thereby establishing a relation between gas flows and gas pressure, which is particularly important for the proper operation of real gas systems.
To evaluate modeling accuracy under the blending pressure problem, we formulated a highly modular and flexible integrated sector-coupled energy system optimization model of the power, natural gas, and (future) hydrogen sectors, where the focus of this paper is on the natural gas and hydrogen sectors and their interconnection with the power sector.
This includes electrolyzer and steam-methane reforming units for hydrogen production, (long-term) natural gas and hydrogen storage units, fuel cell units, co-firing of hydrogen in gas-fired power plants etc. Moreover, the model includes the option to endogenously substitute the natural gas demand of different economic sectors with hydrogen (up to a maximum rate and while ensuring a sufficient energy content for the resulting blend). Thus, the model allows to study integrated expansion planning of generation and transmission infrastructure (power, natural gas and hydrogen) not only in the context of techno-economic characteristics, but also under various policy constraints, e.g. the impact of CO\textsubscript{2} pricing on the deployment of hydrogen in the industry sector.

The presented model was applied in two extensive case studies in an integrated sector-coupled 24-bus IEEE power and 12-node gas system.
In the first case study, we analyzed the impact of the novel blending pressure problem on generation expansion planning in the power and hydrogen sectors, pipeline transmission expansion planning, and operational model results versus a simpler blending transport problem.
Our results indicate that the simple blending transport problem, which omits the physical relation of gas flows and gas pressure in a pipeline, can provide an acceptable approximation for the purpose of generation expansion planning and small hydrogen blending rates. However, expansion planning under the simple framework tends to be overoptimistic, which particularly affects siting decisions not only in the hydrogen but also in the power sector.
Furthermore, omitting the physical realities of gas flows -- particularly in the context of blending -- can strongly impact pipeline transmission expansion planning and the significance of operational results, as the determined gas flows can violate the maximum operating pressure of the gas system to a large extent. Therefore, the planned system might not be able to supply the hydrogen demand.
From these results we conclude that the proposed blending pressure problem, although computationally more intensive than the simpler blending transport problem, is a valuable tool for properly modeling expansion planning in integrated sector-coupled power, natural gas, and hydrogen energy systems.
However, for modeling expansion planning in large scale sector-coupled energy systems, further improving the blending transport problem, e.g. by including constraints to capture pressure characteristics, appears highly relevant.

The second case study centered on the optimal ramp-up of the hydrogen sector on the path towards climate neutrality by leveraging existing pipeline transmission infrastructure for natural gas and hydrogen blending in a totally renewable power system.
Against this background, the optimal deployment of hydrogen is endogenously determined by the presented energy system model across different levels of power and gas sector integration. Our results indicate:
\begin{itemize}
    \item Hydrogen in the power sector: Utilizing the gas transmission system for blending can act as a crucial lever to trigger investments in hydrogen as a power-to-power technology for long-term storage of excess renewable generation -- even for small blending rates. Although the associated reduction of total system cost might be small, this can result in a substantial shift in the cost-optimal investment mix of storage technologies -- in this particular case from battery energy storage systems to hydrogen.
    Despite the investment in hydrogen infrastructure, our results indicate that the number of full-load hours of electrolyzer units is very limited in this case, which does not appear to be economic (from a private investor's perspective).

    \item Hydrogen in the natural gas sector: Exclusively deploying hydrogen as a substitute for natural gas requires sufficiently high costs of natural gas (in addition to the availability of excess renewable power), making this the case with the highest total system costs. However, in contrast to the power-to-power case, the continuity at which hydrogen can be deployed in the gas sector results in a substantially higher number of full-load hours for electrolyzers units and a subordinate role for hydrogen storage.

    \item Sector-coupled power and gas system: In a cost minimization framework, the sector driving the investment in hydrogen infrastructure (in this case study the power sector) is not necessarily the one in which hydrogen is mainly deployed. Furthermore, we observe that the inter-sectoral coupling from the gas towards the power sector is weak with respect to the deployment of hydrogen, while it is more pronounced vice versa. As expected, enhanced sector integration results in the lowest total system costs.
\end{itemize}

In general, the results from our case study strongly indicate that blending could initiate and facilitate the ramp-up of the hydrogen sector and possibly delay the expansion of dedicated hydrogen pipeline infrastructure, e.g. until hydrogen demand is established on a larger scale.
Ultimately, establishing hydrogen demand boils down to hydrogen becoming economic or implementing well-designed policy measures.
For the latter case, our results indicate that the spatio-sectoral distribution of natural gas demand can significantly impact the effectiveness of CO\textsubscript{2} pricing schemes (in terms of total system cost), which could be avoided with other policy schemes, e.g. uniform taxation of natural gas.
The decisive factor for deploying hydrogen in the gas sector, however, is the permissible substitution rate. Against this background, projections on expected hydrogen substitution rates, associated costs for the adaption of the underlying processes, and the correlation with the cost of natural gas are highly relevant topics for future research.

\section*{Acknowledgements}
\label{sec:Acknowledgements}
T. Klatzer gratefully acknowledges the Department of Industrial Economics and Technology Management (I\O T) at the Norwegian University of Science and Technology for hosting his research visit and the funding granted by the Rudolf Chaudoire Foundation and the Erasmus+ program.

\clearpage

\allowdisplaybreaks
\section*{Appendix}
\label{sec:Appendix}
\section*{Nomenclature}

Acronyms:

\begin{longtable}{ l l }
 EU                     & European Union  \\ 
 ESM                    & Energy system model \\
 LP                     & Linear program \\
 MILP                   & Mixed-integer linear program \\
 TP                     & Transport problem \\
 UC                     & Unit commitment \\
 DC-OPF                 & Direct current optimal power flow \\
 SN                     & Single node \\
 LEGO model             & Low-carbon Expansion Generation Optimization model \\
 OF                     & Objective function \\
 OM                     & Operation and Maintenance \\
 CCGT/OCGT              & Combined cycle/Open cycle gas turbine\\
 EL                     & Electrolyzer \\
 SMR                    & Steam-methane reforming \\
 FC                     & Fuel cell \\
 SOC                    & State of charge \\
 MOV                    & Moving window \\
 CU                     & Compressor unit \\
 S-TP                   & Standard transport problem \\
 B-TP                   & Blending transport problem \\
 B-PP                   & Blending pressure problem \\
 INC                    & Incremental \\
 LHS/RHS                & Left-hand/right-hand side \\
 GEP/TEP                & Generation/Transmission expansion planning \\
 H-TEP                  & Hydrogen transmission expansion planning \\
 VRE                    & Variable renewable energy \\
 BESS                   & Battery energy storage system \\
 AEL                    & Alkaline electrolysis \\
 PEMEL                  & Proton exchange membrane electrolyzer \\
 SMR-CCS                & Steam-methane reforming with carbon capture and storage \\
 SOFC                   & Solid oxide fuel cell \\
 CAPEX                  & Capacity expenditure \\
 TRL                    & Technology readiness level \\
 FLH                    & Full-load hour \\
 ETS                    & Emission trading system \\
 NUTS                   & Nomenclature des unités territoriales statistiques \\
 TTF                    & Title Transfer Facility
\end{longtable}

Indices:

\begin{longtable}{ l l }
 $p,pp$                 & Time periods (usually hours)  \\ 
 $rp$                   & Representative periods (usually days)  \\ 
 $k$                    & Time periods within a representative period (usually hours) \\ 
 $\Gamma(p,rp,k)$       & Mapping of periods with representative periods $rp$ and $k$  \\  
 $g$                    & Generating units \\  
 $t(g)$                 & Subset of thermal generation units \\
 $s(g)$                 & Subset of storage generation units \\
 $r(g)$                 & Subset of renewable generation units \\ 
 $h2u$                  & Hydrogen units \\
 $h2g(h2u)$             & Subset of electrolyzer units \\
 $h2p(h2u)$             & Subset of steam-methane reforming units with carbon capture and storage \\
 $h2f(h2u)$             & Subset of fuel cell units \\
 $h2s(h2u)$             & Subset of hydrogen storage units \\
 $h2um(h2u,m)$          & Hydrogen unit $h2u$ connected to gas node $m$ \\
 $ch4u$                 & Natural gas units \\
 $ch4w(ch4u)$           & Subset of natural gas wells \\
 $ch4s(ch4u)$           & Subset of natural gas storage units \\
 $ch4um(ch4u,m)$        & Natural gas unit $ch4u$ connected to gas node $m$ \\  
 $i,j,ii$               & Bus of transmission network \\
 $c$                    & Circuit in transmission network \\
 $ijcc(i,j,c)$          & Candidate transmission line connecting nodes $i$,$j$ with $c$ \\
 $line(i,j)$            & Indicates if a line exists between nodes $i$ and $j$   \\
 $gi(g,i)$              & Generator $g$ connected to node $i$ \\
 $gm(g,m)$              & Generator $g$ connected to gas node $m$ \\
 $m,n$                  & Node of gas transmission system \\
 $l$                    & Pipeline circuit of gas transmission system \\
 $mnl(m,n,l)$           & Pipelines connecting $m$ with $n$ via $l$ \\
 $mnle(m,n,l)$          & Existing pipeline connecting $m$ with $n$ via $l$ \\
 $mnlc(m,n,l)$          & Candidate pipeline connecting $m$ with $n$ via $l$ \\
 $inc$                  & Increment for linearization of pipeline gas flow \\
 $cmp(m,n,l)$           & Compressor unit connecting $m$ with $n$ via $l$ \\
 $cl$                   & Economic class \\
 $sec$                  & Economic sector \\    
 $cls(cl,sec)$          & Relation among economic classes and sector
\end{longtable}

Parameters:

\begin{longtable}{ l l }
 $D^E_{rp,k,i}$                     & Power demand (GW)  \\
 $D^{Gas}_{rp,k,m,cl}$              & Gas demand (MSm\textsuperscript{3}/h)  \\
 $D^{H2}_{rp,k,m,cl}$               & Dedicated hydrogen demand (MSm\textsuperscript{3}/h)  \\
 $H^{H2}, H^{CH4}$                  & Lower heating value of hydrogen and natural gas (GWh/MSm\textsuperscript{3})  \\
 $\underline{SR}^{H2}_{cl}, 
  \overline{SR}^{H2}_{cl}$          & Lower and upper limit for hydrogen substitution (p.u.) \\
 $\underline{B}^{H2},   
  \overline{B}^{H2}$                & Lower and upper limit for hydrogen pipeline blending (p.u.) \\
 $W^{RP}_{rp}$                      & Weight of the representative period (h) \\
 $W^{K}_{k}$                        & Weight of each $k$ within the representative period (h) \\
 $C^{SU}_g$                         & Start-up cost of unit (M\euro) \\
 $C^{UP}_g$                         & Commitment cost of unit (M\euro/h) \\
 $C^{VAR}_g$                        & Variable cost of energy (M\euro/GWh) \\
 $C^{OM}_g$                         & Operation and maintenance cost power unit (M\euro/GWh) \\
 $C^{OM}_{h2u}$                     & Operation and maintenance cost hydrogen unit (p.u.) \\
 $C^{OM}_{ch4u}$                    & Operation and maintenance cost natural gas unit (p.u.) \\ 
 $C^{INV}_g$                        & Investment cost power unit (M\euro/GW/y) \\
 $C^{INV}_{h2g}$                    & Investment cost electrolyzer unit (M\euro/GW/y) \\
 $C^{INV}_{h2p}$                    & Investment cost steam-methane reforming unit (M\euro/(MSm\textsuperscript{3}/h)/y) \\
 $C^{INV}_{h2f}$                    & Investment cost fuel cell unit (M\euro/(MSm\textsuperscript{3}/h)/y) \\
 $C^{INV}_{h2s}$                    & Investment cost hydrogen storage unit (M\euro/(MSm\textsuperscript{3}/h)/y) \\
 $C^{INV}_{ch4u}$                   & Investment cost natural gas unit (M\euro/(MSm\textsuperscript{3}/h)/y) \\
 $C^{L,INV}_{i,j,c}$                & Line investment cost (M\euro/GW/y) \\
 $C^{Pipe,Inv}_{m,n,l}$             & Pipeline investment cost (M\euro) \\
 $C^{CH4}$                          & Cost of natural gas (M\euro/MSm\textsuperscript{3}) \\
 $C^{ENS}$                          & Cost of electricity non-supplied (M\euro/GWh) \\
 $C^{CH4NS}$                        & Cost of natural gas non-supplied (M\euro/MSm\textsuperscript{3}) \\
 $C^{H2NS}$                         & Cost of hydrogen non-supplied (M\euro/MSm\textsuperscript{3}) \\
 $C^{CO2}$                          & Cost of CO\textsubscript{2} allowance (M\euro/MtCO\textsubscript{2}) \\
 $CS^{SU}_g$                        & Start-up gas consumption of unit (GWh) \\
 $CS^{UP}_g$                        & Commitment gas consumption of unit (GWh/h) \\
 $CS^{V}_g$                         & Generation gas consumption of unit (p.u.) \\
 $CS^{Cmp,CH4}_{m,n,l}$             & Compressor natural gas consumption (p.u.) \\
 $CS^{Cmp,H2}_{m,n,l}$              & Compressor hydrogen consumption (p.u.) \\
 $HPE_{h2g}$                        & Hydrogen per unit of electricity (MSm\textsuperscript{3}/GWh) \\
 $HPC_{h2p}$                        & Hydrogen per unit of natural gas (p.u.) \\
 $\overline{P}^E_{g}$               & Technical maximum of power unit (GW) \\
 $\overline{P}^E_{h2g}$             & Technical maximum of electrolyzer unit (GW) \\
 $\overline{P}^{H2}_{h2p}$          & Technical maximum of steam-methane reforming unit (MSm\textsuperscript{3}/h) \\
 $\overline{P}_{ch4u}$              & Technical maximum of natural gas unit (MSm\textsuperscript{3}/h) \\
 $\overline{P}^{H2}_{h2s}$          & Technical maximum production of hydrogen storage unit (MSm\textsuperscript{3}/h) \\
 $\underline{R}^{H2}_{h2s}$         & Technical minimum of hydrogen storage unit (p.u.) \\
 $\overline{CS}^{H2}_{h2s}$         & Technical maximum consumption of hydrogen storage unit (MSm\textsuperscript{3}/h) \\
 $ETP_{h2s}$                        & Energy to power ratio of hydrogen storage unit (hours) \\
 $InRes^{H2}_{h2s,p}$               & Initial reserve of long-term hydrogen storage unit (MSm\textsuperscript{3}) \\
 $\eta^{CH}_{h2s},  
  \eta^{DIS}_{h2s}$                 & Charging and discharging efficiency of hydrogen storage unit (p.u.) \\
 $\overline{F}^{Gas}_{m,n,l}$       & Technical maximum pipeline capacity (MSm\textsuperscript{3}/h) \\
 $F'_{inc,m,n,l}$                   & Function value of linearized gas flow ((MSm\textsuperscript{3}/h)\textsuperscript{2}) \\
 $F_{inc,m,n,l}$                    & Discrete value of linearized gas flow (MSm\textsuperscript{3}/h) \\
 $R^{Gas}_{m,n,l}$                  & Pipeline factor ((MSm\textsuperscript{3}/h)\textsuperscript{2}/bar\textsuperscript{2}) \\
 $\overline{P}^{Sqr}_{m}$           & Technical maximum gas pressure at node (bar\textsuperscript{2}) \\
 $\Lambda^{sqr}_{m,n,l}$            & Compression ratio of compressor unit (p.u.) \\
 $\overline{\Lambda}_{m,n,l}$       & Maximum compression of compressor unit (bar) \\
 $Re$                               & Reynolds number (--) \\
 $D$                                & Nominal pipeline diameter (m) \\
 $v^M$                              & Average gas velocity (m/s) \\
 $\rho^N, \rho^M$                   & Standard and average gas density (kg/Sm\textsuperscript{3}) \\
 $\eta^M$                           & Average dynamic gas viscosity ($1\cdot10^{-6}$ kgs/ms) \\
 $\lambda$                          & Pipeline friction factor (--) \\
 $L$                                & Pipeline length (m) \\
 $\pi$                              & The number $\pi$ (--) \\
 $T^N, T^M$                         & Standard and average gas temperature (K) \\
 $p^N$                              & Standard pressure (bar) \\
 $K^M$                              & Average gas compressibility (--) \\
 $\epsilon$                         & Pipeline roughness (mm) \\
 $\overline{F}^{Cmp,Gas}_{m,n,l}$   & Technical maximum gas flow through compressor (MSm\textsuperscript{3}/h) \\
 $E_g$                              & CO\textsubscript{2} emissions of power unit (MtCO\textsubscript{2}/MSm\textsuperscript{3}) \\
 $E_{cl}$                           & CO\textsubscript{2} emissions of sectoral class (MtCO\textsubscript{2}/MSm\textsuperscript{3}) \\
 $\kappa$                           & Minimum clean production (p.u.) \\
 $EU_g$                             & Indicator of existing power unit (integer) \\
 $EU^{H2}_{h2u}$                    & Indicator of existing hydrogen unit (integer) \\
 $EU^{CH4}_{ch4u}$                  & Indicator of existing natural gas unit (integer) \\
 $\overline{X}_g$                   & Maximum amount of power units to be built (integer) \\
 $\overline{X}^{H2}_{h2u}$          & Maximum amount of hydrogen units to be built (integer) \\
 $\overline{X}^{CH4}_{ch4u}$        & Maximum amount of natural gas units to be built (integer) \\
 $\overline{X}^L_{i,j,c}$           & Maximum amount of transmission lines to be built $\in \{0,1\}$ \\
 $\overline{X}^{Pipe}_{m,n,l}$      & Maximum amount of pipelines to be built $\in \{0,1\}$ \\
 $M$                                & Large positive constant 
\end{longtable}

Variables:

\begin{longtable}{ l l }
 $p^E_{rp,k,g}$             & Power generation of the unit (GW)  \\ 
 $p^{CH4}_{rp,k,ch4w}$      & Natural gas production of the unit (MSm\textsuperscript{3}/h) \\ 
 $p^{H2}_{rp,k,h2u}$        & Hydrogen production of the unit (MSm\textsuperscript{3}/h) \\ 
 $cs^E_{rp,k,g}$            & Power consumption of the power unit (GW)  \\
 $cs^E_{rp,k,h2u}$          & Power consumption of the hydrogen unit (GW)  \\
 $cs^{CH4}_{rp,k,h2u}$      & Natural gas consumption of the hydrogen unit (MSm\textsuperscript{3}/h)  \\
 $cs^{H2}_{rp,k,h2u}$       & Hydrogen consumption of the hydrogen unit (MSm\textsuperscript{3}/h)  \\
 $cs^{CH4,E}_{rp,k,g}$      & Natural gas consumption for power generation of the unit (MSm\textsuperscript{3}/h)  \\
 $cs^{H2,E}_{rp,k,g}$       & Hydrogen consumption for power generation of the unit (MSm\textsuperscript{3}/h)  \\
 $cs^{CH4,Aux}_{rp,k,g}$    & Natural gas consumption for startup and commitment of the unit (MSm\textsuperscript{3}/h)  \\
 $cs^{H2,Aux}_{rp,k,g}$     & Hydrogen consumption for startup and commitment of the unit (MSm\textsuperscript{3}/h)  \\
 $d^{CH4}_{rp,k,m,cl}$      & Variable natural gas demand in the gas sector (MSm\textsuperscript{3}/h)  \\
 $d^{H2}_{rp,k,m,cl}$       & Variable hydrogen demand in the gas sector (MSm\textsuperscript{3}/h)  \\
 $intra^{H2}_{rp,k,h2u}$    & Intra-period state of charge of the hydrogen unit (MSm\textsuperscript{3})  \\
 $inter^{H2}_{p,h2u}$       & Inter-period state of charge of the hydrogen unit (MSm\textsuperscript{3})  \\
 $y_{rp,k,g}$               & Startup decision of the unit (integer)  \\ 
 $u_{rp,k,g}$               & Dispatch commitment of the unit (integer)  \\ 
 $pns_{rp,k,i}$             & Power non-supplied (GW) \\ 
 $h2ns_{rp,k,m,cl}$         & Hydrogen non-supplied (MSm\textsuperscript{3}/h) \\  
 $ch4ns_{rp,k,m,cl}$        & Natural gas non-supplied (MSm\textsuperscript{3}/h) \\ 
 $x_{g}$                    & Investment in power generation capacity (integer)  \\
 $x^{H2}_{h2u}$             & Investment in hydrogen capacity (integer)  \\
 $x^{ch4}_{ch4u}$           & Investment in natural gas capacity (integer)  \\
 $x^{Pipe}_{m,n,l}$         & Investment in pipeline capacity (integer)  \\
 $x^{L}_{i,j,c}$            & Investment in power line capacity (integer)  \\
 $f^{Gas}_{rp,k,m,n,l}$     & Pipeline gas flow (MSm\textsuperscript{3}/h)  \\
 $f^{CH4}_{rp,k,m,n,l}$     & Pipeline natural gas flow (MSm\textsuperscript{3}/h)  \\
 $f^{H2}_{rp,k,m,n,l}$      & Pipeline hydrogen flow (MSm\textsuperscript{3}/h)  \\
 $p^{sqr}_{rp,k,m}$         & Pressure at gas node (bar\textsuperscript{2})  \\
 $\rho_{rp,k,m,n,l}$        & Slack variable (MSm\textsuperscript{3}/h)  \\
 $\delta_{rp,k,inc,m,n,l}$  & Gas flow linking variable (binary)  \\
 $\gamma_{rp,k,inc,m,n,l}$  & Gas flow increment-filling variable (continuous)  \\
 $\alpha_{rp,m,n,l}$        & Pipeline flow direction (binary)  \\
 $f^{Cmp,CH4}_{rp,k,m,n,l}$ & Compressor natural gas flow (MSm\textsuperscript{3}/h)  \\
 $f^{Cmp,H2}_{rp,k,m,n,l}$  & Compressor hydrogen flow (MSm\textsuperscript{3}/h)  \\
\end{longtable}

\bibliography{main}

\end{document}